%% file: chiralII.tex
\input abbrev

\ref\CartanII{ H. Cartan, La Transgression dans un groupe de Lie
et dans un espace fibr\'e principal, Colloque de Topologie,
C.B.R.M., Bruxelles 57-71 (1950).} 
\ref\DKV{M. Duflo, S. Kumar, and M. Vergne, Sur la Cohomologie \'Equivariante des Vari\'et\'es Diff\'erentiables, Ast\'erisque 215 (1993).}
\ref\FS{E. Frenkel and M. Szczesny, Chiral de Rham Complex and Orbifolds, math.AG/0307181.}
\ref\FMS{ D. Friedan, E. Martinec, and S. Shenker,  Conformal Invariance, Supersymmetry and String Theory, Nucl. Phys. B271 (1986) 93-165.} 
 \ref\GMS{V. Gorbounov, F. Malikov, and V. Schectman, Gerbes of Chiral Differential Operators, Math. Res. Lett. 7 (2000) No. 1, 55-66.} 
 \ref\GS{ V. Guillemin and S. Sternberg, Supersymmetry and
Equivariant de Rham Theory, Springer, 1999.}
 \ref\HS{P. Hilton, and U. Stammbach, A Course in Homological
Algebra, 2nd ed., Springer-Verlag, New York 1997.}
\ref\LLI{B. Lian and A. Linshaw, Chiral Equivariant Cohomology I, Adv. Math. 209, 99-161 (2007).}
 \ref\MS{F. Malikov, and V. Schectman, Chiral de Rham Complex, II, math.AG/9901065.}
 \ref\MSV{F. Malikov, V. Schectman, and A. Vaintrob, Chiral de Rham Complex, Comm. Math. Phys, 204, (1999) 439-473.}
\ref\Witten{E. Witten, Two-Dimensional Models With (0,2) Supersymmetry: Perturbative Aspects, Adv. Theor. Math. Phys. 11 (2007) No. 1, 1-63.}

\centerline{\titlefont Chiral Equivariant Cohomology II}

\bs
\centerline{Bong H. Lian, Andrew R. Linshaw and Bailin Song}
\bs

\baselineskip=13pt plus 2pt minus 2pt
ABSTRACT.  This is the second in a series of papers on a new equivariant cohomology that takes values in a vertex algebra. In an earlier paper, the first two authors gave a construction of the cohomology functor on the category of $O(\gs\gg)$ algebras. The new cohomology theory can be viewed as a kind of ``chiralization'' of the classical equivariant cohomology, the latter being defined on the category of $G^*$ algebras a l\`a H. Cartan. In this paper, we further develop the chiral theory by first extending it to allow a much larger class of algebras which we call $\gs\gg[t]$ algebras. In the geometrical setting, our principal example of an $O(\gs\gg)$ algebra is the chiral de Rham complex $\cQ(M)$ of a $G$ manifold $M$. There is an interesting subalgebra of $\cQ(M)$ which does not admit a full $O(\gs\gg)$ algebra structure but retains the structure of an $\gs\gg[t]$ algebra, enough for us to define its chiral equivariant cohomology. The latter then turns out to have many surprising features that allow us to delineate a number of interesting geometric aspects of the $G$ manifold $M$, sometimes in ways that are quite different from the classical theory.

\baselineskip=15pt plus 2pt minus 1pt
\parskip=\baselineskip

\centerline{\bf Contents}\nobreak\medskip{\baselineskip=12pt
 \parskip=0pt\catcode`\@=11 \input toc.tmp \catcode`\@=12 \bigbreak\bigskip}  

\np
\headline{\ifodd\pageno\rightheadline\else\leftheadline\fi}
\def\rightheadline{\tenrm\hfil Chiral Equivariant Cohomology II
\hfil\folio}
\def\leftheadline{\tenrm\folio\hfil B.H. Lian, A.R. Linshaw \& B. Song\hfil}

\newsec{Introduction}

We fix a compact connected Lie group $G$ with a complexified Lie algebra $\gg$.

In section 5 of \LLI, the notion of the chiral equivariant cohomology of an $O(\gs\gg)$ algebra is defined. This is a functor from the category of $O(\gs\gg)$ algebras (where the morphisms are basic differential vertex algebra homomorphisms) to the category of vertex algebras. It turns out that
the functor continues to make sense if we replace the $O(\gs\gg)$ algebra structure with roughly half of that structure. 

Recall that given a Lie algebra $\gg$, we have a Lie superalgebra $\gs\gg=\gg\triangleright\gg^{-1}$ with bracket $[(\xi,\eta),(x,y)]=([\xi,x],[\xi,y]-[x,\eta])$.
See section 3 of \LLI. We then form a loop algebra $\gs\gg[t,t^{-1}]$, equipped with a derivation $d:(\xi,\eta)t^n\mapsto(\eta,0)t^n$, and its associated differential vertex algebra $O(\gs\gg)$. The Lie subalgebra $\gs\gg[t]$ of the loop algebra is clearly preserved by $d$.
Recall that an $O(\gs\gg)$ algebra is a degree-weight graded differential vertex algebra (DVA) $(\cA,d_\cA)$ equipped with a DVA homomorphism $O(\gs\gg)\ra(\cA,d_\cA)$. In particular, $\cA$ is a module over the loop algebra $\gs\gg[t,t^{-1}]$.

What we call an $\gs\gg[t]$ algebra is a vertex algebra equipped with a linear action of the differential Lie superalgebra $\gs\gg[t]$. We notice that the notion of the (chiral) basic subalgebra of an $O(\gs\gg)$ algebra defined in \LLI~ uses only the $\gs\gg[t]$ module structure, rather than the full $O(\gs\gg)$ structure. The module structure is enough for us to define the basic subalgebra of any $\gs\gg[t]$ algebra, hence its basic cohomology. This observation allows us to construct a much more interesting cohomology theory than the one insisting on a full $O(\gs\gg)$ structure. For example, if $M$ is a $G$ manifold, it was shown in \LLI ~that the smooth chiral de Rham complex $\cQ(M)$ is an $O(\gs\gg)$ algebra. It turns out that $\cQ(M)$ contains an interesting subcomplex $\cQ'(M)$ which admits only an $\gs\gg[t]$ structure. Thus our extended theory of chiral equivariant cohomology can now be a potentially useful tool for studying the geometric aspects of the $G$ manifold via the algebra $\cQ'(M)$, in addition to $\cQ(M)$.

Various aspects of the chiral de Rham complex and its generalizations have been studied and developed in recent years. In \MS, the sheaf cohomology was considered in terms of representation theory of affine Lie algebras. In \GMS, cohomological obstructions to the existence of certain subsheaves and global sections were considered, while in \FS, the chiral de Rham complex was studied in the context of finite group actions and orbifolds. 
More recently, chiral differential operators have been interpreted in terms of certain twisted two-dimensional supersymmetric sigma models \Witten. It has been suggested that chiral equivariant cohomology for $G$ manifolds could be a gauged version of a sigma model.$^1$\footnote{}{We thank E. Witten for drawing our attention to this development.} Hopefully, a simple physical interpretation can be found for some of our results on the chiral equivariant cohomology of $\cQ(M)$ and $\cQ'(M)$ (such as the vanishing of the positive weight operators in $\H_G(\cQ(M))$ in some cases.)

{\it Acknowledgement.} B.H.L.'s research is partially supported by a J.S. Guggenheim Fellowship and an NUS grant. A.R.L. and B.S. would like to thank the Department of Mathematics, The National University of Singapore, for its hospitality and financial support during their visit there, where this project was initiated.

\subsec{The smooth ``chiral de Rham sheaf''}

Suppose that $\{\cF_n|~n=0,1,2,\dots\}$ is a family of sheaves of vector spaces on a smooth manifold $M$. In general, the presheaf $\cF$ defined by $\cF(U) = \bigoplus_{n\geq 0} \cF_n(U)$ is not a sheaf. For example, for $M = \R$ and $\cF_n$ a copy of $C^\infty$, if we cover $\R$ by an infinite collection of open intervals, one can use bump functions to construct a family of sections which are compatible on overlaps but do not give rise to a global section of $\cF$. However, $\cF$ does satisfy a slightly weaker version of the reconstruction axiom: \eqn\opencover{
0\ra\cF(U)\rightarrow\prod_i\cF(U_i)\rightrightarrows\prod_{i,j}\cF(U_i\cap U_j)}
is exact for {\it finite} open covers $\{U_i\}$ of an open set $U$. A presheaf which satisfies this weaker exactness condition will be called a {\it weak sheaf}.

There is a mistake in terminology in the construction of the smooth chiral de Rham complex given in Section 2.6 of \LLI. For $M = \R^n$, the sequence $$0\ra\cQ(U)\rightarrow\prod_i\cQ(U_i)\rightrightarrows\prod_{i,j}\cQ(U_i\cap U_j) $$ appearing in the proof of Lemma 2.25 is only exact for finite open covers $\{U_i\}$, so the assignment $U\mapsto \cQ(U)$ should really be called a weak sheaf, instead of a sheaf, on $\R^n$. Note that for any open set $U\subset \R^n$, $\cQ(U)$ is a weight graded vertex algebra, 
so that 
$$\cQ(U) = \bigoplus_{m\geq 0} \cQ(U)[m]$$ 
where $\cQ(U)[m]$ is the subspace of conformal weight $m$. For each $m$, the assignment $U\mapsto \cQ(U)[m]$ does define a sheaf on $\R^n$. Moreover, the statements and proofs of Lemmas 2.25-2.30 of \LLI~are correct if we use the sheaf $\cQ[m]$ instead of $\cQ$. This yields a family of sheaves $\{\cQ_M[m]|~ m\geq 0\}$ on any smooth manifold $M$. Finally, $\cQ_M$ is the weak sheaf of vertex algebras defined by 
$$\cQ_M(U)=\bigoplus_{m\geq 0} \cQ_M[m](U).$$
Alternatively, one can also construct $\cQ_M$ by using $\cQ$ as in Section 2.6 of \LLI, but one must work only with finite open covers.

We should {\it not} work with the sheafification $\tilde{\cQ}_M$ of $\cQ_M$ since (again using bump functions) one can show that $\tilde{\cQ}_M$ is not the direct sum of its weighted pieces, and $\tilde{\cQ}(U)$ is not a vertex algebra in general.

The change of terminology from sheaf to weak sheaf has no practical effect on our theory of chiral equivariant cohomology. The weak sheaf $\cQ_M$ is graded by the sheaves $\cQ_M[m]$. Whenever we need to construct a global section of $\cQ_M$ by gluing together local sections, these sections are always homogeneous of finite weight, so we may work inside $\cQ_M[m]$ for some $m$. The $\gs\gg[t]$-algebra $\cQ'(M)=\cQ_M'(M)$ mentioned above will be the global sections of a weak sheaf $\cQ'_M$ which is graded by subsheaves $\cQ'_M[m]\subset \cQ_M[m]$. 

{\it For simplicity, we will drop the word \lq\lq weak" throughout this paper, and use the word sheaf to denote a weak sheaf when no confusion will arise.}

\subsec{Outline of main results}

We begin with a general formulation of the chiral equivariant cohomology on the category of $\gs\gg[t]$ modules in section 2. 

The functoriality of $\cQ'$ and $\H_G(\cQ'(-))$ is proved in the geometric setting in section 3. The main result here is the isomorphism $\cQ'_{M/G}\cong p_*\cQ'_M$ for principal $G$ bundles $p:M\ra M/G$. 

In section 4, we prove that the universal class given by the Virasoro element $\L$, constructed in \LLI, is in fact a conformal structure on the chiral equivariant cohomology for any semisimple group $G$. This gives a powerful way to decide the vanishing of higher weight classes, since a vertex algebra equipped with a conformal structure is a commutative algebra iff the conformal structure is trivial. We show that for any faithful linear representation of a semisimple group $G$, the chiral equivariant cohomology of the algebra $\cQ(M)$ is purely classical. We also prove that if a $G$ manifold $M$ has a fixed point, then the chiral Chern-Weil map for the algebra $\cQ'(M)$ is injective. In another extreme, we prove that if $G$ acts locally freely, then the chiral equivariant cohomology of $\cQ'(M)$ is purely classical.

When $G$ is abelian, we prove in section 5 that the chiral equivariant cohomology of $M$ is purely classical iff $G$ acts locally freely.

In section 6, we prove a reduction theorem that relates $\H_{G\times T}$ with $\H_G$ when $T$ is abelian. We then apply it to show that for any semisimple group $G$ and a torus $T\subset G$, the chiral equivariant cohomology of $\cQ'(G/T)$ is purely classical.










\subsec{Notation}

The reader is assumed to be familiar with the material in section 5 of reference \LLI. Throughout this paper, we shall adhere to most of the notations introduced in that paper. A partial list of frequently used notations is appended below.

{\it Throughout this paper, a manifold $M$ is always assumed to have a finite open cover consisting of coordinate open sets. Unless stated otherwise, a compact group considered here is assumed connected.}

\settabs\+&\quad  ************\quad\quad &\quad N\cr
\+&\quad  $G,H,T$ &\quad compact connected Lie groups. \cr
\+&\quad  $\gg,\gg^*$ &\quad complexified Lie algebra of $G$ and its dual.  \cr
\+&\quad $\bra\xi,\xi'\ket$&\quad the pairing between $\xi\in\gg$ and $\xi'\in\gg^*$.\cr
\+&\quad  $\gs\gg$ &\quad  the differential Lie superalgebra $\gg\triangleright\gg^{-1}$ with differential $(\xi,\eta)\mapsto(\eta,0)$. \cr
\+&\quad  $\gs\gg[t,t^{-1}]$ &\quad the differential loop algebra of $\gs\gg$.  \cr
\+&\quad  $\gs\gg[t]$ &\quad the differential subalgebra of the loop algebra which is polynomial in $t$.  \cr
\+&\quad  $\gs\gg t^n$ &\quad the subspace of monomials $(\xi,\eta)t^n$ of degree $n$ in $\gs\gg[t,t^{-1}]$. \cr
\+&\quad  $O(\gs\gg)$ &\quad the current algebra of $\gs\gg$ with differential ${\bf d}$.  \cr
\+&\quad  $\cE(\gg),\cS(\gg),\cW(\gg)$ &\quad the semi-infinite exterior, symmetric and Weil algebras of $\gg$. \cr
\+&\quad  $b^\xi, c^{\xi'}$ &\quad the standard generators of the vertex algebra $\cE(\gg)$. \cr
\+&\quad  $\beta^\xi, \gamma^{\xi'}$ &\quad the standard generators of the vertex algebra $\cS(\gg)$. \cr
\+&\quad  $\bra\gamma,c\ket$ &\quad the vertex subalgebra generated by the $\gamma^{\xi'},c^{\xi'}$ in $\cW(\gg)$. \cr
\+&\quad $V=\oplus_{p,n\in\Z,} V^p[n]$&\quad  degree-weight graded vector space whose degree $p$ and weight $n$ \cr
\+&\quad   &\quad subspace is $V^p[n]$. \cr
\+&\quad  $QO(V)$ &\quad the space of Fourier operators $a:V\ra V((z))$.   \cr
\+&\quad  $a(z)=\sum a(p)z^{-p-1}$&\quad power series expansion of a Fourier operator $a$ with $p$th Fourier mode $a(p)$.
  \cr
\+&\quad  $a\circ_p b$ &\quad the $p$th circle product of Fourier operators $a,b$.  \cr
\+&\quad  $\partial a=a\circ_{-2}1$ &\quad the formal derivative of a Fourier operator. \cr

\+&\quad  $\cV\ra QO(\cV)$, $a\mapsto\hat a$&\quad the left regular module of a vertex algebra $\cV$; $\hat a(p)b=a(p)b=a\circ_p b$.  \cr
\+&\quad $(\cA,d)$  &\quad differential vertex algebra (DVA), or an $O(\gs\gg)$ algebra, or an $\gs\gg[t]$ algebra.  \cr
\+&\quad $\cA_{inv},\cA_{hor},\cA_{bas}$ &\quad the chiral invariant, horizontal and basic subalgebras of an $\gs\gg[t]$ algebra $\cA$  \cr
\+&\quad $\Omega=\Omega_M$  &\quad the de Rham sheaf of differential forms of a smooth manifold $M$.  \cr
\+&\quad  $\cQ=\cQ_M$ &\quad the chiral de Rham sheaf of a smooth manifold $M$.  \cr
\+&\quad $d_\cQ$  &\quad the chiral de Rham differential on $\cQ$.  \cr
\+&\quad $b^i,c^i,\beta^i,\gamma^i$ &\quad generators of $\cQ(U)$ associated to a coordinate open set $U$ with\cr
\+&\quad   &\quad   coordinates $\gamma:U\ra\R^n$.  \cr
\+&\quad $a\mapsto a|U$ &\quad  a restriction map $\cF(M)\ra\cF(U)$ of a sheaf $\cF$.  \cr
\+&\quad  $\cQ'=\cQ_M'$ &\quad the sheaf of vertex subalgebras generated by the subsheaf $\Omega=\cQ[0]$ in $\cQ$.  \cr
\+&\quad $\H_{bas}(\cA)$ &\quad the chiral basic cohomology $=H(\cA_{bas},d)$ of an $\gs\gg[t]$ algebra $(\cA,d)$.  \cr
\+&\quad $\H_G(\cA)$ &\quad the chiral equivariant cohomology $=H([\cW(\gg)\otimes\cA]_{bas},d_\cW+d_\cA)$ \cr
\+&\quad  &\quad of an $\gs\gg[t]$ algebra $(\cA,d_\cA)$.  \cr
\+&\quad $f':\cQ_N'\ra f_*\cQ_M'$  &\quad the sheaf homomorphism induced by a smooth map $f:M\ra N$.  \cr
\+&\quad $\H_G(\C){\br\kappa_G\over\ra}\H_G(\cA)$  &\quad the chiral Chern-Weil map of $\cA$  \cr
\+&\quad $\H_{bas}(\cA){\br\pi_G\over\ra}\H_G(\cA)$ &\quad the chiral principal map of $\cA$. \cr
\+&\quad $H_{bas}(A)$  &\quad the classical basic cohomology of a $G^*$ algebra $(A,d)$.   \cr
\+&\quad $H_G(A)$ &\quad the classical equivariant cohomology of a $G^*$ algebra $(A,d)$.  \cr
\+&\quad $W=W(\gg)$ &\quad the classical Weil algebra of $\gg$.  \cr
\+&\quad $(\xi,\eta)\mapsto L_\xi+\iota_\eta$  &\quad an $O(\gs\gg)$ structure $\gs\gg\ra\cA$ on a DVA $(\cA,d)$. \cr
\+&\quad $(\xi,\eta)\mapsto \Theta_\cW^\xi+b^\eta$  &\quad the $O(\gs\gg)$ structure on $\cW(\gg)$; $\Theta_\cW^\xi=\Theta_\cE^\xi+\Theta_\cS^\xi$. \cr
\+&\quad $g$, $L=dg$  &\quad a chiral contracting homotopy and a Virasoro structure on an\cr \+&\quad   &\quad $O(\gs\gg)$ algebra $(\cA,d)$.  \cr
\+&\quad $\xi'\mapsto\theta_{\xi'}$ &\quad chiral connection forms parameterized by $\gg^*\ra\cA[0]$ on an $O(\gs\gg)$ algebra.  \cr
\+&\quad $\xi'\mapsto\Gamma_{\xi'}=g\circ_0\theta_{\xi'}$ &\quad the weight one elements associated to the chiral connection forms.  \cr

\newsec{Chiral Equivariant Cohomology of $\gs\gg[t]$ Modules}

Recall that a differential vertex algebra (DVA) is a degree graded vertex (super)algebra $\cA^*=\bigoplus_{p\in\Z}\cA^p$ equipped with a vertex algebra derivation $d$ of degree 1 such that $d^2=0$.  In particular, for any homogeneous $a,b\in\cA$, 
$$
d(a\circ_n b)=(da)\circ_n b+(-1)^{deg~a}a\circ_n(db).
$$
As in \LLI, {\it the DVAs considered in this paper will have an additional $\Z$-valued weight gradation, which is always assumed to be bounded below by 0 and compatible with the degree, i.e. $\cA^p=\bigoplus_{n\geq0}\cA^p[n]$.}

\definition{An $\gs\gg[t]$ module is a degree-weight graded complex $(\cA^*,d_\cA)$ equipped with a Lie algebra homomorphism $\rho:\gs\gg[t]\ra End~\cA^*$, such that  for all $x\in\gs\gg[t]$ we have
\item{$\bullet$} $\rho(dx)=[d_\cA,\rho(x)]$;
\item{$\bullet$} $\rho(x)$ has degree 0 whenever $x$ is even in $\gs\gg[t]$, and degree -1 whenever $x$ is odd, and has weight $-n$ if $x\in\gs\gg t^n$.}

\definition{Given an $\gs\gg[t]$ module $(\cA,d)$, we define the chiral horizontal, invariant and basic subspaces of $\cA$ to be respectively
$$\eqalign{
\cA_{hor}&=\{a\in\cA|\rho(x)a=0~\forall x\in\gg^{-1}[t]\},\cr
\cA_{inv}&=\{a\in\cA|\rho(x)a=0~\forall x\in\gg[t]\},\cr
\cA_{bas}&=\cA_{hor}\cap\cA_{inv}.
}$$
An $\gs\gg[t]$ module $(\cA,d)$ is called an $\gs\gg[t]$ algebra if it is also a DVA such that $\cA_{hor},\cA_{inv}$ are both vertex subalgebras of $\cA$.}

When we wish to emphasize the role of the group $G$, we shall use the notations $\cA_{G-hor}$, $\cA_{G-inv}$, $\cA_{G-bas}$, and we use the term $G$ chiral instead of chiral.

Given an $\gs\gg[t]$ module $(\cA,d)$, it is clear that $\cA_{inv}$, $\cA_{bas}$ are both subcomplexes of $\cA$. But $\cA_{hor}$ is not a subcomplex of $\cA$ in general.
Note also that the Lie algebra $\gs\gg[t]$ is not required to act by derivations on a DVA $\cA$ to make it an $\gs\gg[t]$ algebra.

If $(\cA,d)$ is an $O(\gs\gg)$ algebra and $(\xi,\eta)\mapsto L_\xi+\iota_\eta$ denotes the $O(\gs\gg)$ structure, then any subDVA $\cB$ which is closed under the operators $(L_\xi+\iota_\eta)\circ_p$, $p\geq0$, is an $\gs\gg[t]$ algebra.


\others{Example}{An interesting subsheaf of $\cQ$.} Let $M$ be a $G$ manifold. Recall that in \LLI, the smooth chiral de Rham complex $\cQ_M$ on $M$ is a direct sum of sheaves $\cQ_M[m]$, constructed by transferring sheaves $\cQ[m]$ on $\R^n$ with the ``Grothendieck topology.'' We define a subsheaf $\cQ'[m]$ of $\cQ[m]$ on $\R^n$ by
$$
\cQ'[m](U)=\bra \Omega(U)\ket[m],
$$
i.e. the weight $m$ subspace of the vertex subalgebra in $\cQ(U)$ generated by $\Omega(U)=\cQ[0](U)$. Then it is easy to check that $\cQ'(U):=\bigoplus_{m\geq0}\cQ'[m](U)$ is an abelian vertex algebra which is linearly isomorphic to $\C[\partial^i\gamma^j,\partial^i c^j|i>0,j=1, ... ,n]\otimes\Omega(U)$. Then by the same construction as in Section 2.6 of \LLI, we can transfer the (weak) sheaf $\cQ'$ on $\R^n$ to $M$ and obtain a (weak) sheaf of vertex algebras on $M$. Thus for each point $x\in M$, we have
$$
\cQ_M'[m]_x=\bra(\Omega_M)_x\ket[m]\subset\cQ_M[m]_x.
$$

An $O(\gs\gX)$ structure is defined on $\cQ(M)$ (see Remark 3.3 of \LLI). We observe that the space of sections of $\cQ'$ is closed under the chiral de Rham differential $d_\cQ$, which means that $\cQ'$ is a sheaf of DVAs. More importantly, for any vector field $X$, the operators $L_X\circ_p$ and $\iota_X\circ_p$ preserve $\cQ'$ for all $p\geq0$. Since the $G$ action on $M$ induces a homomorphism $O(\gs\gg)\ra O(\gs\gX)$, it follows that $\cQ'(M)$ is also closed under the action of the Lie algebra $\gs\gg[t]$, and $\cQ'(M)_{inv}=\cQ(M)_{inv}\cap\cQ'(M)$, $\cQ'(M)_{hor}=\cQ(M)_{hor}\cap\cQ'(M)$ are both subalgebras of $\cQ'(M)$. Thus $(\cQ'(M),d_\cQ)$ is an $\gs\gg[t]$ algebra. Note that by Lemma 2.9 of \LLI, the $\gs\gg[t]$ action on $\cQ'(M)$ is in fact by derivations because $\cQ'(M)$ is abelian.

Note that the assignment $M\mapsto\cQ(M)$ is {\it not} functorial in the category of manifolds, because its construction involves both differential forms and vector fields on $M$ at the same time. But we shall see in the next section that $M\mapsto\cQ'(M)$ is, in fact, functorial. Roughly speaking, the main reason is that $M\mapsto\Omega(M)$ is functorial, and that $\cQ'(M)$ is a vertex algebra generated canonically by the commutative algebra $\Omega(M)$.

Let $\cA,\cB$ be $\gs\gg[t]$ algebras. A DVA homomorphism $f:\cA\ra\cB$ is said to be {\it chiral basic} if $f(\cA_{bas})\subset\cB_{bas}$. Likewise for a module homomorphism. If $f:\cA\ra\cB$ is chiral basic and $\cC$ is an $O(\gs\gg)$ algebra, then it need not be true that the homomorphism $1\otimes f:\cC\otimes\cA\ra\cC\otimes\cB$ is chiral basic. The case we will study is when $\cC=\cW(\gg)$. The most important case is when $\cA=\C$, for which this is not a problem.
For if $\cB,\cC$ are any $\gs\gg[t]$ algebras, the canonical DVA homomorphism $\cC\ra\cC\otimes\cB$, $c\mapsto c\otimes 1$, is chiral basic because
$$
\cC_{bas}\otimes1\subset(\cC\otimes\cB)_{bas}.
$$

The following generalizes the notion of the chiral basic cohomology of an $O(\gs\gg)$ algebra. See Definition 5.1 of \LLI.

\definition{Let $(\cA^*,d)$ be an $\gs\gg[t]$ module. Its chiral basic cohomology $\H_{bas}(\cA)$ is the cohomology of the complex $(\cA_{bas}^*,d)$.}

\definition{For any $\gs\gg[t]$ module $(\cA^*,d_\cA)$, we define its chiral equivariant cohomology $\H_{G}(\cA)$ to be the chiral basic cohomology of the tensor product
$$
(\cW(\gg)\otimes\cA,d_\cW+d_\cA)
$$ 
where $d_\cW=J(0)\otimes 1+K(0)\otimes1$ and $d_\cA\equiv 1\otimes d_\cA$. We call the canonical maps 
$$
\kappa_G=\kappa^\cA_G:\H_{bas}(\cW(\gg))=\H_G(\C)\ra\H_G(\cA),~~~
\pi_G=\pi_G^\cA:\H_{bas}(\cA)\ra\H_G(\cA),
$$ 
induced by $\cW(\gg)\ra\cW(\gg)\otimes\cA$ and $\cA\ra\cW(\gg)\otimes\cA$, respectively, the chiral Chern-Weil map and the chiral principal map of $\cA$.} 

In this paper, the two main examples of chiral equivariant cohomology $\H_G(\cA)$ are given by the $\gs\gg[t]$ algebras $\cA=\cQ(M),\cQ'(M)$ of a $G$ manifold $M$. We will be especially interested in their chiral Chern-Weil and principal maps.
The name for the second map is motivated by the following classical example. If $M$ is a $G$ manifold, then $\Omega(M)$ is a $G^*$ algebra. The canonical $G^*$ algebra homomorphism $\Omega(M)\ra W\otimes\Omega(M)$ induces a map $H_{bas}(\Omega(M))\ra H_G(M)$. Note that $(\Omega(M)_{bas},d_{dR})$ is a subcomplex of the de Rham complex, hence all positive degree elements in $H_{bas}(\Omega(M))$ are nilpotent. Now suppose $M$ is a principal $G$ bundle. Then the map $H_{bas}(\Omega(M))\ra H_{bas}(W\otimes\Omega(M))=H_G(M)$ is an isomorphism, and both sides coincide with $H(M/G)$. See section 4.3 of \GS. In general, the map measures how far $M$ is from being a principal bundle.


The $\gs\gg[t]$ algebra $\cA=\cQ'(M)$ has a special feature. The equipped $\gs\gg[t]$ structure actually comes from an $O(\gs\gg)$ structure on a larger DVA, namely $\cB=\cQ(M)$, of which $\cA$ is a subalgebra. Because of this special feature, we have
$$
\cA_{bas}=\cA^{\gs\gg[t]}=\cB_{bas}\cap\cA.
$$
Since $\cB_{bas}$ is a commutant subalgebra of $\cB$, $\cA_{bas}$ is a subalgebra of $\cA$.

\subsec{The small chiral Weil model for $\gs\gt[t]$ algebras}

In this subsection, we assume that $G=T$ is abelian.

In \LLI, we introduced for any $O(\gs\gt)$ algebras $\cA$, the small chiral Cartan model
$$
\cC_T(\cA)=\bra\gamma\ket\otimes\cA_{inv}
$$
with differential $d_T=1\otimes d_\cA-(\gamma^{\xi_i}\otimes\iota_{\xi_i})(0)$. Here $\bra\gamma\ket$ is the vertex subalgebra generated by the $\gamma^{\xi'}$ in $\cW(\gt)$. We have shown that $H(\cC_T(\cA),d_T)\cong\H_T(\cA)$. In this subsection, we will describe the analogue of this in the chiral Weil model and for $\gs\gt[t]$ algebras $\cA$ by using the chiral Mathai-Quillen isomorphism $\Phi$. {\it We shall assume that the $\gs\gg[t]$ structure of $\cA$ is inherited from an $O(\gs\gg)$ algebra $\cA'$ containing $\cA$.} As usual, we denote the $O(\gs\gg)$ structure by $(\xi,\eta)\mapsto L_\xi+\iota_\eta$. Let $\bra\gamma,c\ket$ be the vertex subalgebra generated by the $\gamma^{\xi'},c^{\xi'}$ in $\cW(\gt)$.

\lemma{$\cC_T(\cA)=\Phi(\bra\gamma,c\ket\otimes\cA_{inv})_{hor}$.}
\proof
Recall that $\Phi^{-1}=e^{-\phi(0)}$ where
\eqn\dumb{
\phi(0)=\sum c^{\xi_i}(-p-1)\otimes\iota_{\xi_i}(p).
}
Since the $c^{\xi}(-p-1)$ kill $\bra\gamma\ket$ for $p<0$, we only need to sum over $p\geq0$ above, if we let $\Phi^{-1}$ act on $\cC_T(\cA)$. For abelian $G=T$, the $L_\xi, \iota_\eta$ commute so that each $\iota_{\xi}(p)$ preserves $\cA_{inv}$. This shows that $\Phi^{-1}\cC_T(\cA)\subset \bra\gamma,c\ket\otimes\cA_{inv}$. Since elements of $\Phi^{-1}\cC_T(\cA)\subset(\cW(\gt)\otimes\cA')_{bas}$ are chiral basic, and since $O(\gs\gt)$ acts trivially on $\cW(\gt)$ because $G=T$ is abelian, it follows that
\eqn\dumbI{
\Phi^{-1}\cC_T(\cA)\subset(\bra\gamma,c\ket\otimes\cA_{inv})_{bas}=(\bra\gamma,c\ket\otimes\cA_{inv})_{hor}.
}

By a property of the chiral Mathai-Quillen map, $$\Phi(\bra\gamma,c\ket\otimes\cA_{inv})_{hor}\subset \Phi(\cW\otimes\cA_{inv})_{bas}\subset\cW_{hor}\otimes\cA_{inv}.$$
But since the sum \dumb~ is over $p\geq0$, we have $\Phi(\bra\gamma,c\ket\otimes\cA_{inv})_{hor}\subset \bra\gamma,c\ket\otimes\cA_{inv}$. Since $\cW_{hor}\cap\bra\gamma,c\ket=\bra\gamma\ket$, it follows that
$$
\Phi(\bra\gamma,c\ket\otimes\cA_{inv})_{hor}\subset\bra\gamma\ket\otimes\cA_{inv}=\cC_T(\cA).
$$
Together with \dumbI, this implies our assertion. $\Box$

This shows that
$$
\cD_T(\cA){\br def\over =}\Phi^{-1}\cC_T(\cA)=(\bra\gamma,c\ket\otimes\cA_{inv})_{hor}
$$
with the differential 
$$
\Phi^{-1}d_T\Phi=d_\cW+d_\cA=(\gamma^{\xi_i}b^{\xi_i})(0)+d_\cA
$$
is isomorphic to the {\it small chiral Cartan model} $(\cC_T(\cA),d_T)$, where $$d_T=-(\gamma^{\xi_i'}\otimes\iota_{\xi_i})(0)+d_\cA.$$ We call $(\cD_T(\cA),d_\cW+d_\cA)$ the {\it small chiral Weil model} for $\cA$. To summarize, we have

\theorem{The maps $\cC_T(\cA){\br\Phi^{-1}\over\ra}\cD_T(\cA)\hra(\cW(\gt)\otimes\cA)_{bas}$ induce the isomorphisms $\H_T(\cA)\cong H(\cD_T(\cA),d_\cW+d_\cA){\br\Phi\over\cong} H(\cC_T(\cA),d_T)$. }
\thmlab\SmallWeilModel

By using the relation $d_\cW\beta^\xi=b^\xi$, we find that the complex $(\bra\beta,b\ket,d_\cW)$ is acyclic. Using the fact that the $\beta^\xi,b^\xi$ are chiral horizontal, one can show directly that the inclusion
$(\bra\gamma,c\ket\otimes\cA_{inv})_{hor}\hra(\cW(\gt)\otimes\cA)_{bas}$ is a quasi-isomorphism of DVAs, giving another proof that the small chiral Weil model for $\cA$ also computes $\H_T(\cA)$.

The preceding theorem also shows that chiral equivariant cohomology is functorial with respect to abelian groups. Let $T,H$ be tori and $\cA$ be an $\gs\gh[t]$ algebra. Let $T\ra H$ be a map of Lie groups, so that $\cA$ has an induced $\gs\gt[t]$ algebra structure. The map $T\ra H$ also induces a map $\bra\gamma,c\ket_H\ra\bra\gamma,c\ket_T$ between the subalgebras generated by the $\gamma,c$ in $\cW(\gh)$ and $\cW(\gt)$, respectively. This gives rise to a DVA map $\cD_H(\cA)\ra\cD_{T}(\cA)$. Thus in cohomology we have a natural vertex algebra homomorphism $\H_H(\cA)\ra\H_T(\cA)$. It is easy to check that this natural map respects compositions.

\newsec{The Object $\cQ'$}


\subsec{Global sections}

\lemma{Each $\cQ'[n]$ is naturally a $C^\infty$-module over $M$, as a sheaf of vector spaces.}
\proof
We have $\cQ'[n]_x=\bra \Omega_x\ket[n]$, which is clearly a module over $C^\infty_x\subset\Omega_x$ with respect to the commutative associative product $\circ_{-1}$. $\Box$ 

The preceding lemma is false for the sheaf $\cQ$. The reason is that if $f,g$ are functions, then $:f(:g\beta:):\neq:(fg)\beta:$ in general.

\lemma{$\cQ'$ is locally free. In fact, if $U$ is a coordinate open set, then $\cQ'(U)$ is generated by $\Omega(U)$ as a vertex algebra. In particular, each weighted piece of the sheaf $\cQ'$ is a vector bundle of finite rank.}
\proof
Obviously $\bra\Omega(U)\ket\subset\cQ'(U)$. From the construction of the MSV sheaf $\cQ$, it is clear that each $a\in\bra\Omega_x\ket=\cQ'_x$ can be represented in local coordinates near $x$ as an element of $\C[\partial^{k+1}\gamma^i,\partial^kc^i|i=1, ... ,n;~k=0,1,... ]\otimes C^\infty(U)=\bra\Omega(U)\ket$. Thus we have the reverse inclusion $\cQ'(U)\subset\bra\Omega(U)\ket$. 

Finally, for a given weight, there are only a finite number of weight homogeneous monomials in the variables $\partial^{k+1}\gamma^i,\partial^kc^i$, and this number is the same on every small open set $U$. This shows that each weighted piece of the sheaf $\cQ'$ is locally free of finite rank. $\Box$

Let $\cF$ be a sheaf of vector spaces on $M$, $U$ an open subset in $M$, and $a\in\cF(U)$. Let $V$ be the largest open subset $V\subset U$ such that $a|V=0$. The support of $a$ is defined as the closure of $U\backslash V$ in $M$ and is denoted by $supp~a$. If $supp~a\subset U$, then there is a unique extension of $a$ to a global section $\tilde a\in\cF(M)$ with the same support as $a$. 
{\it We call $\tilde a$ the global extension of $a$ by zero.} If $\cK\subset\cF(U)$ is the subspace of elements $a$ with $supp~a\subset U$, then the image of the extension map $\cK\ra\cF(M)$, $a\mapsto\tilde a$, is precisely the subspace in $\cF(M)$ consisting of elements with support in $U$.

Suppose $\cF$ is a $C^\infty$-module.  A family of global sections $a_\alpha\in\cF(M)$ is said to be locally finite if for every $x\in M$, there is a neighborhood $U_x\ni x$ such that $a_\alpha|U_x=0$ for all but finitely many $\alpha$. In this case $\sum_\alpha a_\alpha$ is well-defined. It is the unique global section $a$ such that $a|U_x=\sum_\alpha (a_\alpha|U_x)$. It is independent of the choice of the $U_x$. Let $U_\alpha$ be open sets forming a locally finite cover of $M$. Let $\rho_\alpha$ be a partition of unity subordinate to that cover. Then given $a\in\cF(M)$, the family $\rho_\alpha a$ is locally finite, and $\sum_\alpha\rho_\alpha a=a$.

Assuming, still, that $\cF$ is a $C^\infty$-module on $M$, then given an open set $U$, if $a\in\cF(U)$ and if $\rho\in C^\infty(M)$ with $supp~\rho\subset U$, then $(\rho|U)a\in\cF(U)$ has support in $U$. We denote by $\rho a$ its global extension by zero.

\lemma{Given an open set $U$, let $\cK$ be the subspace of $\cQ'(U)$ consisting of elements with support in $U$. Then $\cK$ is an ideal in $\cQ'(U)$, and extension by zero $\cK\ra\cQ'(M)$ preserves all circle products and $\partial$.}
\proof
Since the restriction maps for the sheaf $\cQ'$ are vertex algebra homomorphisms, the ideal property of $\cK_U$ is clear. Let $\cJ$ be the subspace of $\cQ'(M)$ consisting of elements with support in $U$. Now the extension map $\cK\ra\cQ'(U)$ is a linear isomorphism onto $\cJ$, with inverse given by the restriction map, i.e. $a\mapsto\tilde a\mapsto\tilde a|U=a$. In particular for $a,b\in\cK$, $\widetilde{a\circ_p b}|U=a\circ_p b$. But since the restriction maps of the sheaf $\cQ'$ are vertex algebra homomorphisms, we also have $(\tilde a\circ_p\tilde b)|U=a\circ_p b$. Since $\widetilde{a\circ_p b}$ and $\tilde a\circ_p\tilde b$ both lie in $\cJ$, and they have the same image under the restriction $\cJ\ra\cK$, they must be equal. This shows that extension by zero preserves all circle products. 

Note that $\partial a=a\circ_{-2}1$, and the proof that extension by zero preserves $\partial$ is similar. $\Box$


\lemma{Let $U\subset M$ be a coordinate open set with coordinates $\gamma^i$ and $c^i=d\gamma^i$. If $f\in C^\infty(U)$ has $supp~f\subset U$, then each $f\partial^k\gamma^i, f\partial^k c^i$ lies in $\bra\Omega(M)\ket$.}
\proof
We will use induction to show that $f\partial^k\gamma^i$ lies in $\cA=\bra\Omega(M)\ket$; for $f\partial^k c^i$ the proof is similar. The case $k=0$ is trivial. Suppose $f\partial^k\gamma^i\in\cA$ for all $f\in C^\infty(U)$ with $supp~f\subset U$. Note that all the derivatives of such an $f$ have supports in $U$. We have that
$$
\partial(f\partial^k\gamma^i)=\partial f\partial^k\gamma^i+f\partial^{k+1}\gamma^i
$$
lies in $\cA$. Choose a function $g:U\ra\R$ with $g|F=1$ and $supp~g\subset U$, where $F=supp~\partial f\subset U$. Since the extensions of $\partial f, g\partial^k\gamma^i$ lie in $\cA$, and $\partial f\partial^k\gamma^i=\partial f ~g\partial^k\gamma^i$, it follows that the extension of $f\partial^{k+1}\gamma^i$ also lies in $\cA$. $\Box$

\theorem{$\Omega(M)$ generates $\cQ'(M)$ as a vertex algebra.}
\proof\thmlab\GlobalSectionQPrime
The idea is this: if $\Omega(M)$ is to generate $\cQ'(M)$, it had better generate all those elements in $\cQ'(M)$ obtained by extending local sections by zero using bump functions. On the other hand, every global section is a sum of these extensions, by partition of unity. Now the proof.

One inclusion $\bra\Omega(M)\ket\subset\cQ'(M)$ is obvious. Let $a\in\cQ'(M)$ be a homogeneous element of weight $m$. Let $\{U_\alpha\}$ be a locally finite cover of coordinate open sets of $M$. Let $\rho_\alpha$ be the partition of unity with $supp~\rho_\alpha\subset U_\alpha$. We have $a=\sum_\alpha\rho_\alpha a$. Since $\cQ'[m]$ is a sheaf, this possibly infinite sum makes sense. It suffices to show that each $\rho_\alpha a\in\bra \Omega(M)\ket$.

Note that $supp~\rho_\alpha a\subset U_\alpha$. Thus it suffices to show that for any coordinate open set $U$ and any vertex operator $a\in\cQ'(M)$ with $supp~a\subset U$, we have $a\in\bra\Omega(M)\ket$. We can write $a|U$ as a sum of monomials in local coordinates with coefficients in $C^\infty(U)$. That $supp~a\subset U$ means that all the coefficient functions appearing in $a|U$ have supports in $U$. Since $\widetilde{a|U}=a$, to show that $a\in\bra\Omega(M)\ket$, it suffices to show that for each monomial $A$ in $\C[\partial^{k+1}\gamma^i,\partial^kc^i|i=1,..,n;~k=0,1,...]$, and each function $\rho\in C^\infty(U)$ with $F=supp~\rho\subset U$, we have $\rho A\in\bra\Omega(M)\ket$. 

Let $f: U\ra\R$ be a smooth function with $f|F=1$ and $supp~f\subset U$. Then $f\rho=\rho$.
Thus if we modify the monomial $A$ by replacing each $\partial^{k+1}\gamma^i,\partial^kc^i$ by $f\partial^{k+1}\gamma^i,f\partial^kc^i$, respectively, we get a vertex operator $B\in\bra\Omega(M)\ket$ such that $\rho B=\rho A$. But each extension $f\partial^{k+1}\gamma^i,
f\partial^k c^i$ lies in $\bra\Omega(M)\ket$, by the preceding lemma. It follows that $\rho A=\rho B\in\bra\Omega(M)\ket$. This completes the proof. $\Box$

\subsec{Functoriality of $\cQ'$}

\theorem{The assignment $M\mapsto\cQ_M'$ is functorial in the sense that any smooth map $\varphi:M\ra N$ induces a homomorphism of sheaves of DVAs, $\varphi':\cQ_N'\ra\varphi_*\cQ_M'$ over $N$. Moreover, if $M{\br\varphi\over\ra} N{\br\psi\over\ra} P$ are smooth maps, then the map $\cQ_P'{\br(\psi\circ\varphi)'\over\ra}(\psi\circ\varphi)_*\cQ_M'$ 
coincides with the composition of
$\cQ_P'{\br\psi'\over\ra}\psi_*\cQ_N'{\br\psi_*\varphi'\over\ra}\psi_*\varphi_*\cQ_M'$. }

\proof\thmlab\FunctorialityQPrime
When no confusion arises, we shall drop the subscript from $\cQ'_M$. It will be clear that when $U$ is an open set of $M$, the section space $\cQ'(U)$ will, of course, refer to that of the sheaf $\cQ_M'$.

Recall that for any coordinate open set $U$ of $M$, $\cQ'(U)$ is the abelian vertex subalgebra generated by $\Omega(U)$ in $\cQ(U)$. Let $\varphi:M\ra N$ be a smooth map of manifolds. Let $U,V$ be two coordinate open sets in $M,N$, respectively. If $U\subset\varphi^{-1} V$ then we define $\varphi':\cQ'(V)\ra\cQ'(U)$ by extending the classical pullback $\varphi^*:\Omega(V)\ra\Omega(U)$ as follows. 
By construction of $\cQ(U)$ (see section 2.4 of \LLI), we find that the relations in $\cQ'(U)$ are precisely
$$
1(z)=id,~~~(ab)(z)=~:a(z)b(z):~,~~\forall a,b\in\Omega(U),~~~\partial g={\partial g\over\partial\gamma^i}~\partial\gamma^i,~~\forall g\in C^\infty(U),
$$
where $\gamma:U\ra\R^m$ is any coordinate chart on $U$. 
Likewise for the algebra $\cQ'(V)$. This shows that there is a unique extension $\varphi':\cQ'(V)\ra\cQ'(U)$ such that $\varphi'|\Omega(V)=\varphi^*$, and $\partial[\varphi'(y)]=\varphi'(\partial y)$ for all $y\in\cQ'(V)$. Note that we have the convention $\cQ'(\emptyset)=0$.

The differential $d_\cQ:\cQ'\ra\cQ'$ is characterized by the property that it is a derivation of the circle products (in particular commutes with $\partial$) such that $d_\cQ|\Omega(U)$ is the de Rham differential. Since $\varphi^*$ commutes with de Rham and $\varphi'$ commutes with $\partial$, it follows that $\varphi'$ commutes with $d_\cQ$.

If $U_1,U_2\subset\varphi^{-1}V$ are two coordinate open sets in $M$, then it is easy to check that the induced maps above are compatible with restrictions in $M$, i.e. they fit into the commutative diagram
 \eqn\dumbII{\matrix{
& & \cQ'(U_1)\cr
&\nearrow& \downarrow\cr
\cQ'(V) & \ra & \cQ'(U_1\cap U_2)\cr
&\searrow& \uparrow\cr
& &\cQ'(U_2)\cr
}}
the vertical arrows being restriction maps.
This shows that the induced maps $\cQ'(V)\ra\cQ'(U)$, where $U$ ranges over all coordinate open subsets of $\varphi^{-1}V$, patch together to form a map $\varphi':\cQ'(V)\ra\cQ'(\varphi^{-1}V)=(\varphi_*\cQ')(V)$.
Likewise, if $V_1,V_2\supset\varphi(U)$ are coordinate open sets in $N$, then the induced maps are compatible with restrictions in $N$, and they fit into the commutative diagram
\eqn\dumbI{\matrix{
\cQ'(V_1) & &\cr
\downarrow& \searrow&\cr
\cQ'(V_1\cap V_2)  & \ra & \cQ'(U)\cr
\uparrow& \nearrow&\cr
 \cQ'(V_2)& &\cr
}}
This shows that we have a well-defined homomorphism $\varphi':\cQ_N'\ra\varphi_*\cQ_M'$ of sheaves of DVAs over $N$.

To check that the induced homomorphisms $\varphi',\psi'$ respect compositions, it suffices to do it locally. That is, if $U,V,W$ are coordinate open sets in $M,N,P$, respectively such that
$\varphi(U)\subset V$, $\psi(V)\subset W$, then the composition of the maps 
$$
\cQ'(W)\ra\cQ'(V)\ra\cQ'(U)
$$ 
induced locally by $\psi,\varphi$, respectively, agrees with the map $\cQ'(W)\ra\cQ'(U)$ induced locally by $\psi\circ\varphi$. But this follows immediately from three facts: (1) the classical pullbacks 
$$
\Omega_P{\br(\psi\circ\varphi)^*\over\ra}\psi_*\varphi_*\Omega_M,~~~
\Omega_P{\br\psi^*\over\ra}\psi_*\Omega_N,~~~
\Omega_N{\br\varphi^*\over\ra}\varphi_*\Omega_M
$$
have this composition property; (2) the induced maps $\varphi',\psi'$ are vertex algebra homomorphisms; (3) the local sections spaces $\cQ'(U),\cQ'(V),\cQ'(W)$ are vertex algebras generated by the subspaces $\Omega(U),\Omega(V),\Omega(W)$, respectively. The details are tedious but straightforward and are left to the reader. $\Box$

Let $M,N$ be $G$ manifolds. For $h\in G$, let $h_M,h_N$ be the corresponding diffeomorphisms on $M,N$, respectively. Let $\varphi:M\ra N$ be a $G$ equivariant map. Then
$$
\varphi\circ h_M=h_N\circ\varphi:M\ra N.
$$
By the preceding theorem, we have
$$
(\varphi_*h_M')\circ\varphi'=(\varphi\circ h_M)'=(h_N\circ\varphi)'=
({h_N}_*\varphi')\circ h_N'.
$$
Hence we have a commutative diagram of sheaves over $N$
$$\matrix{
\varphi_*\cQ_M' &\br\varphi_* h_M'\over\lra & {h_N}_*\varphi_*\cQ_M'=\varphi_* {h_M}_*\cQ_M'\cr
\varphi'\uparrow & & \uparrow{h_N}_*\varphi'\cr
\cQ_N' & \br h_N'\over\lra & {h_N}_*\cQ_N'
}$$
Thus on global sections over $N$, we have the commutative diagram
$$\matrix{
\cQ_M'(M) &\br h_M'\over\lra & \cQ_M'(M)\cr
\varphi'\uparrow & & \uparrow\varphi'\cr
\cQ_N'(N) & \br h_N'\over\lra & \cQ_N'(N)
}$$
which shows that the DVA homomorphism $\varphi'$ is $G$ equivariant. In the next section, we shall see that this map is in fact a homomorphism of $\gs\gg[t]$ algebras.



\subsec{Chiral horizontal and basic conditions are local}

\lemma{Let $\cP$ be a subsheaf of a sheaf of vertex algebras $\cV$, and $\tau$ a global section of $\cV$. For every open set $U$, let $\cV_\tau(U)=Com(\tau|U,\cP(U))$. Then $\cV_\tau$ defines a subsheaf of $\cV$.}
\proof
The space $\cV_\tau(U)$ consists of elements $a\in\cP(U)$ which commute with $\tau|U$. Thus it is a vertex subalgebra of $\cP(U)$. It is enough to show that this commutant condition is local, i.e. two vertex operators $a,b\in\cV(U)$ commute iff for every $x\in U$, there is a neighborhood $V\ni x$ in $U$ such that $a|V,b|V$ commute. Recall that two vertex operators $a,b$ commute iff $a\circ_p b=0$ for all $p\geq0$. Since restrictions to open subsets are vertex algebra homomorphisms, i.e. they preserve circle products, it follows that $a,b$ commute iff $a|V,a|V$ commute for every open subset $V\subset U$. This shows that the commutant condition is local. $\Box$

\corollary{Let $\cP$ be a subsheaf of vertex algebras of $\cQ=\cQ_M$ which is stable under the $\gs\gg[t]$ action. Then $U\mapsto\cP(U)_{hor}$ and $U\mapsto\cP(U)_{bas}$ define two subsheaves of vertex algebras of $\cP$.}
\proof
The chiral horizontal and chiral basic conditions are both commutant conditions defined by
global sections $L_\xi,\iota_\xi\in\cQ(M)$, $\xi\in\gg$. $\Box$

\subsec{Principle $G$ Bundles}

Throughout this subsection, let $M$ be a principal $G$ bundle. We denote by
$$
p:M\ra M/G
$$
the natural projection. Think of $p$ as a $G$ equivariant map, where $G$ acts trivially on $M/G$. So we have a sheaf homomorphism of $\gs\gg[t]$ algebras $p':\cQ'_{M/G}\ra p_*\cQ'_M$, by the functoriality of $\cQ'$. 

\lemma{$p'\cQ'_{M/G}\subset (p_*\cQ'_M)_{bas}=p_*(\cQ'_M)_{bas}$.}
\proof
Since $p'$ is a chiral basic homomorphism, we have $p'(\cQ'_{M/G})_{bas}\subset (p_*\cQ'_M)_{bas}$. But we have $(\cQ'_{M/G})_{bas}=\cQ'_{M/G}$ because $G$ acts trivially on $M/G$, and $$(\cQ'_M)_{bas}(p^{-1}U)=\cQ'_M(p^{-1}U)_{bas}$$ for any open set $U$ in $M/G$ because the chiral basic condition is local. This completes the proof. $\Box$

\theorem{$p':\cQ'_{M/G}\ra p_*(\cQ'_M)_{bas}$ is a sheaf isomorphism.}
\thmlab\SheafGReduction

Restricting this isomorphism to weight zero, we recover the well-known fact that $\Omega_{M/G}\cong p_*(\Omega_M)_{bas}$. Thus this theorem is a vertex algebra analogue of the classical theorem.

Since this is a local question, it suffices to show that for every small open set $U\subset M/G$, $p':\cQ'(U)\ra\cQ'(p^{-1}U)_{bas}$ is an isomorphism.
For small $U\subset M/G$, we have $p^{-1}U\cong U\times G$. Under this identification, $p:U\times G\ra U$ is just the projection map. So we can regard $p':\cQ'(U)\ra\cQ'(p^{-1}U)_{bas}$ as $p':\cQ'(U)\ra\cQ'(U\times G)_{bas}$, which is induced by the projection map $U\times G\ra U$. Let $j:U\ra U\times G$, $x\mapsto (x,e)$. This is not $G$ equivariant. But it induces $j':\cQ'(U\times G)\ra\cQ'(U)$, which we can restrict to the basic subalgebra $\cQ'(U\times G)_{bas}$. Note that $p\circ j=id_U$. By functoriality, $j'\circ p'=id_{\cQ'(U)}$. So we have

\lemma{Locally $p'$ is injective.}

By Theorem \GlobalSectionQPrime, for a small open set $U\subset M/G$ ($G$ compact), we have
$$
\cQ'(U)=\bra c^i,C^\infty(U)\ket,
~~~~\cQ'(U\times G)=\bra c^i,\theta_{\xi'},C^\infty(U\times G)\ket
$$ 
where $c^i$ are the coordinate 1-forms on $U$, and the $\theta_{\xi'}$ are the pullbacks of the connection forms on $G$. Since the Fourier modes $\theta_{\xi'}(-p-1),\iota_\xi(p)$, $p\geq0$, generate a copy of Clifford algebra acting on the space $\cQ'(U\times G)$, it follows that
$$
\cQ'(U\times G)_{hor}=\bra c^i, C^\infty(U\times G)\ket.
$$
More generally, for any open set $V\subset G$, we have
\eqn\dumb{
\cQ'(U\times V)_{hor}=\bra c^i, C^\infty(U\times V)\ket.
}

Let's consider the case when $G=T$ is abelian first.

\lemma{$p':\cQ'(U)\ra\cQ'(U\times T)_{bas}$ is surjective.}
\proof
Let $\gamma^i$, $i=1,\dots,n$, be coordinates on $U$. Put $T=\R^m/\Z^m$. For $j=n+1,\dots,n+m$, let $\gamma^j$ be the flat coordinates of $T$, and let $\xi_j$ be the standard basis of $\gt=\R^m$. Then $X_{\xi_j}={\partial\over\partial\gamma_j}$ which are global vector fields. It follows that $L_{\xi_j}=\beta^j$. Note that while the $\gamma^j$ are 
local coordinates on $T$, the $\partial^k\gamma^j$, $k\geq1$, are all global sections because the transition functions for our coordinates are $\gamma^j\mapsto\gamma^j+1$, as in the case of a circle. This shows that
$$
\cQ'(U\times T)_{hor}=\C[\partial^k c^i,\partial^{k+1} \gamma^i, \partial^{k+1}\gamma^j|i=1,..,n;~j=n+1,..,n+m;~k=0,1,...]\otimes C^\infty(U\times T).
$$
By the Heisenberg relations $\beta^j(z)\gamma^l(w)\sim\delta_{jl}(z-w)^{-1}$, it follows that
$$
\cQ'(U\times T)_{bas}=\C[\partial^k c^i, \partial^{k+1} \gamma^i|i=1,..,n;~k=0,1,...]\otimes C^\infty(U)
$$
where $C^\infty(U)$ consists of functions pulled back from $U$ by the projection map $U\times T\ra U$.
It is now clear that the map $p'$ in this case is surjective. $\Box$

We now consider the surjectivity question of $p'$ in the nonabelian case. For $\xi\in\gg$, the corresponding vector field $X_\xi$ on $G$ is locally of the form $X_\xi=f_{\xi j}{\partial\over\partial\gamma^j}$, where $\gamma^j$ are local coordinates on $G$. If $\xi_i$ form a basis of $\gg$, then the coefficient functions $f_{\xi_i j}$ form an invertible matrix. We denote the inverse matrix by $g_{j\xi_i}$.

\lemma{Define the local operators $\tilde\beta^j=:g_{j\xi_i}L_{\xi_i}:$. Then for $h\in C^\infty(U\times G)$,
$$
\tilde\beta^j(z)h(w)\sim{\partial h\over\partial\gamma^j}~(z-w)^{-1}.
$$
In particular, for a coordinate open set $V\subset G$, we have $\tilde\beta^j\circ_p=\beta^j\circ_p$, for $p\geq0$, acting on $\cQ'(U\times V)_{hor}$. Here the $\beta^j$, $j=1,..,m$, correspond to the coordinate functions $\gamma^j$ on $V$. }
\proof
The OPE calculation is straightforward. The second assertion follows from \dumb. $\Box$

\lemma{Let $a\in\cQ'(U\times G)_{hor}$. Then $\tilde\beta^j\circ_p a=0$ for all $j$ and $p\geq0$ holds on every coordinate open set $U\times V$ iff  $a\in\cQ'(U)$, i.e. $a\in Im~p'$.}
\proof
Applying the preceding lemma locally on $U\times V$ with $h=\gamma^j$, we see that
$\tilde\beta^j\circ_p$ acts by ${\partial\over\partial\gamma^j(-p-1)}$ on functions. From the description of the horizontal subalgebra given above, it follows that those horizontal elements $a$ satisfying $\tilde\beta^j\circ_p a=0$ for all $j$ and $p\geq0$ are exactly those with no dependence on all $\gamma^j(-p-1)$ locally. But those are exactly the elements of $\cQ'(U)$. $\Box$

{\it Proof of Theorem \SheafGReduction:}
We have reduced it to a local question and have shown that $p':\cQ'(U)\ra\cQ'(U\times G)_{bas}$ is injective for small open sets $U\subset M/G$. To prove that this is surjective, suppose $a\in\cQ'(U\times G)_{hor}$ is a chiral basic element, i.e. $L_\xi\circ_p a=0$, $p\geq0$, i.e. $L_\xi$ commutes with $a$ in $\cQ(U\times G)$. In particular, this holds on every coordinate open set $U\times V$. 
Since $a$ also commutes with functions, it commutes with the Wick products $:g_{k\xi_i}L_{\xi_i}:~=\tilde\beta^k$ in $\cQ(U\times V)$. By the preceding lemma, $a|(U\times V)$ lies in the image of $p'$. This completes the proof. $\Box$

\newsec{The Chiral Characteristic Class $\kappa_G(\L)$}


Given an $\gs\gg[t]$ algebra $(\cA^*,d_{\cA})$, we wish to find a condition under which the chiral Chern-Weil map $\kappa_G:\H_G(\C)\ra\H_G(\cA)$ kills all positive weight elements. Likewise for the chiral principal map $\pi_G:\H_{bas}(\cA)\ra\H_G(\cA)$. Then, we apply it to the case $\cA=\cQ'(M)$ for a $G$ manifold $M$. We shall first establish the functoriality of $\H_G$ for $G$ manifolds.

\subsec{Functoriality of $\H_G(\cQ')$}

\lemma{Let $\cA,\cB$ be subalgebras of DVAs $\cA',\cB'$, respectively. Suppose $\cA'$ is equipped with an $O(\gs\gg)$ structure $(\xi,\eta)\mapsto L^\cA_\xi+\iota^\cA_\eta$, and likewise for $\cB'$. If $\iota_\xi^\cA\circ_n\cA\subset\cA$ and $\iota_\xi^\cB\circ_n\cB\subset\cB$ for all $n\geq0$, $\xi\in\gg$, and if $f:\cA\ra \cB$ is a DVA homomorphism such that $f(\iota_\xi^\cA\circ_na)=\iota_\xi^\cB\circ_nf(a)$, then $\cA,\cB$ are $\gs\gg[t]$ algebras such that $f$ and $id_\cW\otimes f$ are both chiral basic homomorphisms.}
\proof\thmlab\BiggerAlgebras
Let $d_\cA$ denote the differential on $\cA'$. Then $d_\cA\iota_\xi^\cA=L_\xi^\cA$, so that $L_\xi^\cA\circ_n\cA\subset\cA$. It follows that $\cA$ inherits an $\gs\gg[t]$ structure from the $O(\gs\gg)$ structure of $\cA'\supset\cA$. 
Likewise for $\cB$. Applying $d_\cB$, the differential on $\cB'$, to
\eqn\dumb{
f(\iota_\xi^\cA\circ_na)=\iota_\xi^\cB\circ_nf(a),
} 
and using the identity $d_\cB f(x)=f( d_\cA x)$, we find that  $f(L_\xi^\cA\circ_na)=L_\xi^\cB\circ_nf(a)$. This shows that $f$ is an $\gs\gg[t]$ algebra homomorphism. It follows that 
$$
f(\cA_{bas})=f(\cA^{\gs\gg[t]})\subset\cB^{\gs\gg[t]}=\cB_{bas},
$$ 
implying that $f$ is chiral basic.

From the identities \dumb, it follows that 
$$
(id_\cW\otimes f)[(b^\xi\otimes 1+1\otimes\iota_\xi^\cA)\circ_nx]
=(b^\xi\otimes 1+1\otimes\iota_\xi^\cB)\circ_n[(id_\cW\otimes f)(x)].
$$
Now $\cW\otimes\cA'\supset\cW\otimes\cA$ has an $O(\gs\gg)$ structure defined by $(\xi,\eta)\mapsto (\Theta_\cW^\xi\otimes 1+1\otimes L_\xi^\cA)+(b^\xi\otimes 1+1\otimes\iota_\xi^\cA)$; likewise for $\cW\otimes\cB'$. So we can apply the result above to the map $id_\cW\otimes f:\cW\otimes\cA\ra\cW\otimes\cB$, and see that this is again a chiral basic homomorphism of $\gs\gg[t]$ algebras. $\Box$


\lemma{Let $a,b,c$ be pairwise circle commuting homogeneous vertex operators. Then for $m,n\in\Z$, we have
$$
a\circ_m(b\circ_nc)=(-1)^{(deg~a)(deg~b)}b\circ_n(a\circ_m c)+\sum_{p\geq0}\left(\matrix{m\cr p}\right)(a\circ_pb)\circ_{m+n-p}c.
$$}
\proof\thmlab\Associativity
Let $\cA$ be the vertex algebra generated by $a,b,c$. Let $\cA\ra QO(\cA)$, $u\mapsto\hat u$, denote the left regular module action. Then by Lemma 2.6 of \LLI,
$$
[\hat a(m),\hat b(n)]=\sum_{p\geq0}\left(\matrix{m\cr p}\right)(\widehat{a\circ_pb})(m+n-p).
$$
Applying both sides to $c$ yields the desired identity. $\Box$


\lemma{Suppose that $f:\cA\ra\cB$ is a DVA homomorphism, and that $S\subset\cA$ generates $\cA$ as a vertex algebra. Let $u,v$ be vertex operators which circle commute with (but do not necessarily lie in) $\cA$ and $\cB$, respectively, such that
\eqn\dumb{
u\circ_n a\in\cA,~~~f(u\circ_n a)=v\circ_n f(a),~~~\forall n\geq0,
}
for $a\in S$. Then \dumb~holds for all $a\in\cA$.}
\proof
By induction, it suffices to show that 
if \dumb~ holds for $a=x$ and $y$, then it holds for $a=x\circ_p y$. This follows easily from the preceding lemma. $\Box$

Let $\varphi:M\ra N$ be an equivariant map of $G$ manifolds. We wish to study the induced DVA homomorphism $\varphi':\cQ'(N)\ra\cQ'(M)$ given by Theorem \FunctorialityQPrime. Denote by $(\xi,\eta)\mapsto L_\xi^M+\iota_\eta^M$ and $\mapsto L_\xi^N+\iota_\eta^N$ the $O(\gs\gg)$ structures on $\cQ(M),\cQ(N)$ respectively \LLI. Recall that $L_\xi^M=d_\cQ\iota^M_\xi$, $L_\xi^N=d_\cQ\iota_\xi^N$. By construction, $\varphi'$ restricted to classical forms $\Omega(N)\ra\Omega(M)$ is just the pullback map, and the operator $\iota_\xi^M\circ_0$ on forms is just the interior multiplication by the vector field $x^M_\xi$ generated by $\xi\in\gg$ on $M$. Likewise for $N$. In particular, we can characterize the $\gs\gg[t]$ structure locally on the subalgebra $\cQ'(U)\subset\cQ(U)$, generated by $\Omega(U)$, by the property
$$
\iota^M_\xi\circ_n\omega=\delta_{n,0}~\iota^M_\xi\omega,~~n\geq0,~\omega\in\Omega(U).
$$
Note that since $\iota_\xi^M\in\cQ(U)$ is a vertex operator, we have 
$$
\partial(\iota_\xi\circ_n a)=-n~\iota_\xi\circ_{n-1} a+\iota_\xi\circ_n\partial a
$$
for all $a\in\cQ'(U)$. 

If $\alpha\in\Omega^{k+1}(N)$, then for $x_1,..,x_k\in T_mM$, we have
$$\eqalign{
[\varphi'(\iota_\xi^N\circ_0\alpha)](x_1,..,x_k)
&=\alpha(x_\xi^N,\varphi_*x_1,..,\varphi_*x_k)\cr
&=\alpha(\varphi_*x_\xi^M,\varphi_*x_1,..,\varphi_*x_k),
~~~x_\xi^N=\varphi_*x_\xi^M~at~\varphi(m)\in N\cr
&=[\iota_\xi^M\circ_0\varphi'(\alpha)](x_1,..,x_k).
}$$
It follows that $\varphi'(\iota_\xi^N\circ_0\alpha)=\iota_\xi^M\circ_0\varphi'(\alpha)$. Since $\alpha\in\Omega(N)$ has weight zero, this equation holds if we replace $\circ_0$ by $\circ_n$, for any $n\geq0$, because both sides would then be zero. Applying the preceding lemma to the case
$f=\varphi'$, $\cA=\cQ'(N)$, $\cB=\cQ'(M)$, $u=\iota_\xi^M$, $v=\iota_\xi^N$, $S=\Omega(N)$, we conclude that
$$
\iota_\xi^N\circ_n\alpha\in\cQ'(N),~~~\varphi'(\iota_\xi^N\circ_n\alpha)=\iota_\xi^M\circ_n\varphi'(\alpha),~~~\forall n\geq0,
$$
for all $\alpha\in\cQ'(N)$, $\xi\in\gg$. Now apply Lemma \BiggerAlgebras~to the case $f=\varphi'$, $\cA=\cQ'(N)$, $\cA'=\cQ(N)$, $\cB=\cQ'(M)$, $\cB'=\cQ(M)$. We have

\theorem{$\varphi':\cQ'(N)\ra\cQ'(M)$ and $1\otimes\varphi':\cW(\gg)\otimes\cQ'(N)\ra\cW(\gg)\otimes\cQ'(M)$ are chiral basic homomorphisms of $\gs\gg[t]$ algebras. In particular, the assignment $M\mapsto\H_G(\cQ'(M))$ is functorial on the category of $G$ manifolds.}

\subsec{Vanishing theorem and half $O(\gs\gg)$ TVA}

Recall that for each $c\in\C$, there exists a vertex algebra $\cV(c)$ with a single generator $L$ such that
\eqn\dumb{
L(z)L(w)\sim\half c(z-w)^{-4}+2L(w)(z-w)^{-2}+\partial L(w)(z-w)^{-1}
}
and that it has the following universal property: if $\cB$ is a vertex algebra containing an element $L'$ with the same OPE, then there exists a (unique) vertex algebra homomorphism $\cV(c)\ra\cB$ such that $L\mapsto L'$. We shall usually denote $L'$ simply as $L$.

\definition{Let $\cA$ be a weight graded vertex algebra. We say that $\cA$ has an (outer) Virasoro structure of central charge $c$ if there exists a vertex algebra $\cA'\supset\cA$ and an element $L\in\cA'$ satisfying \dumb, such that $L\circ_p \cA\subset\cA$ for all $p\geq0$. If, furthermore, $L\in\cA$, then we call $L$ an inner Virasoro structure of $\cA$. If, in addition to this, we have $L\circ_0=\partial$ and $L\circ_1|\cA[n]=n$, then we call $L$ a conformal structure of $\cA$.}

\lemma{Let $L$ be an inner Virasoro structure of a vertex algebra $\cA'$ and $\cA\subset\cA'$ is a subalgebra generated by a set $S$. If 
\eqn\dumb{
L\circ_p a\subset S,~~~p\geq0,
}
for all $a\in S$, then $L$ is a Virasoro structure of $\cA$.}
\proof
We want to show that \dumb~holds for all $a\in\cA$. It suffices to show that if it holds for $a=x,y$, then it holds for all $a=x\circ_p y$. Again, this follows from Lemma \Associativity. $\Box$

\others{Example}{$\beta\gamma$ system revisited.} Recall that (Examples 2.15 and 2.19 of \LLI) we have a vertex algebra $\cA'$ with two even generators $\beta,\gamma$ and relations
$$
\beta(z)\gamma(w)\sim(z-w)^{-1},~~~\beta(z)\beta(w)\sim0,~~~\gamma(z)\gamma(w)\sim0.
$$
It has a conformal structure given by $L=:\beta\partial\gamma:$ with central charge 1, which is not contained in the subalgebra $\cA=\bra\gamma\ket$ generated by $\gamma$. But since $L(z)\gamma(w)\sim\partial\gamma(w)(z-w)^{-1}$, $L$ is a Virasoro structure of $\cA$ which is not inner.

\theorem{Let $(\cA',d)$ be a DVA equipped with an $O(\gs\gg)$ structure, $(\xi,\eta)\mapsto L_\xi+\iota_\eta$, and let $\cA$ be an $\gs\gg[t]$ subalgebra of $\cA'$.
\item{1.} If $\cA$ has no negative weight elements, then $\cA[0]$ is a $G^*$ algebra with product $\circ_{-1}$, and $\H_{bas}(\cA)[0]=H_{bas}(\cA[0])$. 
\item{2.} If $\cA'$ has an inner Virasoro structure $L$ such that the $L_\xi,\iota_\xi$ are primary of weight 1 and if $L\circ_p\cA\subset\cA$ for $p\geq0$, then $L\circ_p\cA_{bas}\subset\cA_{bas}$ for $p\geq0$.
\item{3.} If, furthermore, $L$ is a conformal structure of $\cA'$ which is $d$ closed, then $L\circ_p$ operates on $\H_{bas}(\cA)$ for $p\geq0$, and we have $L\circ_0=\partial$, $L\circ_1=n$ on $\H_{bas}(\cA)[n]$.
\item{4.} If, in addition to the assumptions above, there exists a chiral horizontal element $g\in\cA'$ 
such that $L=d g$ and $g\circ_p\cA\subset\cA$ for $p\geq0$, then $g\circ_p\cA_{bas}\subset\cA_{bas}$ for $p\geq0$, and we have $\H_{bas}(\cA)[n]=0$ for $n>0$.
\item{}}
\proof\thmlab\VirasoroVanish
For part 1, that $\cA[0]$ is a (super) commutative associative algebra with product $\circ_{-1}$ follows from Lemma 2.9 of \LLI. Moreover, $L_\xi\circ_n,\iota_\xi\circ_n$ kill $\cA[0]$ for all $n>0$. Since the $L_\xi\circ_0,\iota_\xi\circ_0$ act on $\cA[0]$ by derivations, they define a $G^*$ structure on $\cA[0]$ such that the classical basic complex $\cA[0]_{bas}$ coincides with the weight zero subpace of the chiral basic complex $\cA_{bas}$. In particular, we have $\H_{bas}(\cA)[0]=H_{bas}(\cA[0])$. 

Part 2 is an immediate consequence of Lemma 5.19 of \LLI. 

Part 3 follows from part 2 and the fact that $L$ is a conformal structure of $\cA'$. 

Consider part 4. Applying $d$ to $\iota_\xi\circ_p g=0$, $p\geq0$, we get
$$
L_\xi\circ_p g=\iota_\xi\circ_p dg=\iota_\xi\circ_p L
=\left\{\matrix{0&p\neq1\cr-\iota_\xi&p=1.}\right.
$$
The last equality holds because $\iota_\xi$ is primary of weight 1 with respect to $L$. This shows that $g\circ_q$ preserves $\cA_{bas}$ for $q\geq0$. If $a\in\cA_{bas}[n]$ is $d$ closed, then $L\circ_1 a=na=d (g\circ_1a)$. Since $g\circ_1a$ is basic, this shows that if $n>0$, then $a$ is $d$ exact in $\cA_{bas}$. This proves part 4.
$\Box$

Note that the obstruction for $g\circ_n$ above to act on cohomology $\H_{bas}(\cA)$ is precisely $L\circ_n$ because $d g=L$. But if the $L\circ_n$ act trivially, then the theorem shows that there would be nothing of positive weight in cohomology, and the theory is essentially classical.

Many of our results are about the class of algebras satisfying the assumptions in Theorem \VirasoroVanish, so it makes sense to name these algebras.



\definition{An $O(\gs\gg)$ topological vertex algebra (TVA) is a differential vertex algebra  $(\cA,d)$ equipped with an $O(\gs\gg)$ structure $(\xi,\eta)\mapsto L_\xi+\iota_\eta$, a chiral horizontal homogeneous element $g$, such that $L=d g$ is a conformal structure, with respect to which the $\iota_\eta$ are primary of weight one. We call $g$ a chiral contracting homotopy of $\cA$. Given an $O(\gs\gg)$ TVA $(\cA,d)$, a differential vertex subalgebra $\cB$ is called a half $O(\gs\gg)$ TVA if the nonnegative Fourier modes of the vertex operators $\iota_\xi$ and $g$ preserve $\cB$. We shall always assume that the weights in $\cA$ are nonnegative.}

Note that in the preceding definition of a half $O(\gs\gg)$ TVA $\cB$, the nonnegative Fourier modes of $L_\xi=d\iota_\xi$ and $L=d g$ automatically preserve $\cB$. We shall often abuse the language, and call $g$ and $L$, respectively, a chiral contracting homotopy and a conformal structure of $\cB$, even though neither element necessarily lies in $\cB$.

\lemma{Let $\cA$ be an $O(\gs\gg)$ TVA with chiral contracting homotopy $g$. Then $g$ and $L=dg$ are $G$ invariant.}
\proof\thmlab\gLGinvariance
Since the $\iota_\xi$ are primary of weight one, so are the $L_\xi$. It follows that the Fourier modes $L_\xi(0),\iota_\xi(0)$ commute with $L$. In particular, $L$ is $G$ invariant. Since a chiral contracting homotopy $g$ is assumed to be chiral horizontal, we have
$$
0=[d,[\iota_\xi(0),g]]=[L_\xi(0),g]-[\iota_\xi(0),L]=[L_\xi(0),g].
$$
This shows that $g$ is $G$ invariant as well. $\Box$


\others{Example}{$\cQ(M)$ and $\cQ'(M)$.}
The algebra $(\cQ(M),d_\cQ)$ is a DVA equipped with a conformal structure $L^\cQ$ and a chiral contracting homotopy $g^\cQ$ such that $d_\cQ g^\cQ=L^\cQ$. When $M$ is a $G$ manifold, then $\cQ(M)$ becomes an $O(\gs\gg)$ TVA because $g^\cQ$ is chiral horizontal. In this case, the subalgebras $\cQ(M),\cQ'(M)\subset\cQ(M)$ are half $O(\gs\gg)$ TVAs. Since $\cQ(M)[0]=\cQ'(M)[0]=\Omega(M)$, Theorem \VirasoroVanish~ yields
$$
\H_{bas}(\cQ(M))=\H_{bas}(\cQ'(M))=H_{bas}(\Omega(M)).
$$
That is, the chiral basic cohomology of both algebras $\cQ(M),\cQ'(M)$ coincide with the classical  basic cohomology of $\Omega(M)$.

\others{Example}{Semi-infinite Weil algebra $\cW(\gg)$.}
This algebra has a conformal structure $L^\cW$ (denoted by $\omega_\cW$ in Lemma 5.4 of \LLI), and an element $g^\cW$ (denoted by $h$ in Lemma 5.13 of \LLI) such that $L^\cW=d_\cW g^\cW$. But $g^\cW$ fails to be chiral horizontal because $b^\xi(z) g^\cW(w)\sim\beta^\xi(w)(z-w)^{-2}$. Thus it satisfies all the assumptions up to parts 1-3 of Theorem \VirasoroVanish, but not all of 4, hence it fails to be a (half) $O(\gs\gg)$ TVA. In fact, $\H_{bas}(\cW)=\H_G(\C)$ always contains nonzero elements of positive weights, as shown in \LLI. This is precisely what makes the chiral Chern-Weil map $\kappa_G:\H_G(\C)\ra\H_G(\cA)$ of an $\gs\gg[t]$ algebra $\cA$ interesting!

\lemma{Let $\cA$ be a vertex algebra and $a\in\cA$. If $\partial a=0$, then $a$ commutes with all of $\cA$.}
\proof\thmlab\PartialZero
Let $b\in\cA$. To see that $a,b$ commute, it suffices to show that in the left regular module we have $\hat a\circ_p\hat b=0$ for $p\geq0$. Now the image of this vertex operator under the creation map $\hat\cA\ra\cA$, $\hat c\mapsto c$, is $\hat a(p)b$. But $\partial\hat a=0$ means that the Fourier modes $\hat a(p)=0$ for all $p\neq -1$. In particular $\hat a(p)b=0$ for $p\geq0$, hence $\hat a\circ_p\hat b=0$. $\Box$

\corollary{If $\cA$ is a half $O(\gs\gg)$ TVA, and if a class $x\in\H_G(\cA)$
is in the image of the chiral principal map $\pi_G:\H_{bas}(\cA)\ra\H_G(\cA)$, then $\partial x=0$. In fact, any class $x\in\H_G(\cA)$ of the form $x=[1\otimes a]$ of positive weight is zero. Thus $\pi_G\H_{bas}(\cA)$ commutes with all of $\H_G(\cA)$ in this case.}
\proof
Such classes are exactly those in the image of the chiral principal map. But $\H_{bas}(\cA)$ has no positive weight elements by Theorem \VirasoroVanish. The last assertion follows from the preceding lemma. $\Box$

\corollary{If $M$ is a principal $G$ bundle, then the vertex subalgebra generated by the weight zero elements in $\H_G(\cQ(M))$ or $\H_G(\cQ'(M))$ is the commutative algebra $H(M/G)$ with no positive weights.}
\proof
For $\cA=\cQ(M)$ or $\cQ'(M)$, the algebra $\cW\otimes\cA$ has no negative weight elements, so that by Theorem \VirasoroVanish,
$$
\H_G(\cA)[0]=H_G(\cA[0])=H_G(\Omega(M))=H(M/G).
$$
On the other hand this coincides with $\pi_G H_{bas}(\Omega(M))=\pi_G\H_{bas}(\cA)$,
which is of course a vertex algebra. This shows that the vertex algebra generated by $\H_G(\cA)[0]$ in $\H_G(\cA)$ is just $H(M/G)$. $\Box$

Later, this result will be refined to show that for arbitrary $G$ whose action on $M$ is locally free, $\H_G(\cQ(M)),\H_G(\cQ'(M))$ have no positive weight elements at all.



\theorem{Suppose that $\cA$ is a half $O(\gs\gg)$ TVA and that $G$ is simple. Then $\H_G(\cA)=H_G(\cA[0])$ iff $\kappa_G(\L)=0$.}
\proof\thmlab\VirasoroVanishII
By Theorem \VirasoroVanish, part 1,
$$
\H_G(\cA)[0]=\H_{bas}(\cW\otimes\cA)[0]=H_{bas}(W\otimes\cA[0])=H_G(\cA[0]).
$$
In particular, since $\kappa_G(\L)$ has weight two, the only-if part of our assertion follows. 
Now suppose $\L\otimes 1\in(\cW\otimes\cA)_{bas}$ represents the zero class in $\H_G(\cA)$, i.e. it is equal to $-dx$ for some $x\in(\cW\otimes\cA)_{bas}$ where $d=d_\cW+d_\cA$. Recall that (see the proof of Theorem 7.16 of \LLI) $\L=d_\cW(\Theta_\cS^{\xi_i}b^{\xi_i}+\beta^{\xi_i}\partial c^{\xi_i})$. This implies that 
$$
a=(\Theta_\cS^{\xi_i}b^{\xi_i}+\beta^{\xi_i}\partial c^{\xi_i})\otimes1+x
$$
is $d$ closed. Since $x$ is basic and $\Theta_\cS^{\xi_i}b^{\xi_i}\otimes1$ is obviously horizontal, i.e. killed by $(b^\xi\otimes 1+1\otimes\iota_\xi)\circ_n$ for all $n\geq0$, it follows that $\beta^{\xi_i}\partial c^{\xi_i}\otimes1-a$ is horizontal. So,
$$
h=\beta^{\xi_i}\partial c^{\xi_i}\otimes1-a+1\otimes g^\cA
$$
is also horizontal. In particular $h\circ_1$ acts in $\cW\otimes\cA$ preserving horizontal elements. We have $dh=d_\cW(\beta^{\xi_i}\partial c^{\xi_i})\otimes1+1\otimes d_\cA g^\cA=L^\cW\otimes1+1\otimes L^\cA$. It follows that  
$$
[d,h\circ_1]|(\cW\otimes\cA)[n]=-n~id.
$$
Thus if $y$ is any $d$ closed basic element of weight $n$, we have
$$
d(h\circ_1 y)=-n y.
$$
For $n>0$, this says that $d z=y$ for some horizontal element $z$. But this implies that $z$ is also basic.
Hence $y$ represents the zero class in $\H_G(\cA)$. $\Box$

\theorem{If $\cA$ is a half $O(\gs\gg)$ TVA and if $G$ is semisimple, then
$$
\kappa_G(\L)\circ_n=(L^\cW\otimes1+1\otimes L^\cA)\circ_n,~~~n\geq0.
$$ 
In particular, $\kappa_G(\L)$ is a conformal structure on the vertex algebra $\H_G(\cA)$.}
\proof\thmlab\ConformalStructureTheorem
Apply Theorem \VirasoroVanish~with the structure $\cW\otimes\cA\subset\cW\otimes\cA'\ni$
$L^\cW\otimes1+1\otimes L^\cA$, playing the respective roles of $\cA\subset\cA'\ni$ $L$ there.
Then $(L^\cW\otimes1+1\otimes L^\cA)\circ_n$, $n\geq0$, operate on $\H_G(\cA)$.
The second assertion now follows from the first assertion.

Put $L=L^\cW\otimes1+1\otimes L^\cA$, $d=d_\cW+d_\cA$. In the proof of Theorem 7.16 of \LLI~(which considers the case of simple $G$, but which generalizes readily to semisimple $G$ with minor changes), we have
$$
\L=d_\cW(\Theta_\cS^{\xi_i}b^{\xi_i}+\beta^{\xi_i}\partial c^{\xi_i})=
d_\cW(\Theta_\cS^{\xi_i}b^{\xi_i})+L^\cW.
$$
Since $L^\cA=d_\cA g^\cA$, this shows that
$$
\kappa_G(\L)-L=dh,~~~h=\Theta_\cS^{\xi_i}b^{\xi_i}\otimes 1-1\otimes g^\cA.
$$
Note that $h$ is chiral horizontal. Suppose $a\in(\cW\otimes\cA)_{bas}$ is $d$ closed. Then for $n\geq0$
$$
[\kappa_G(\L)-L]\circ_n a=(dh)\circ_n a=d(h\circ_n a).
$$
Since the left side is chiral horizontal, and $h\circ_n a$ is horizontal, it follows that $h\circ_n a$ is also basic. So the left side is $d$ exact in $(\cW\otimes\cA)_{bas}$. $\Box$

Taking $\cA=\C$, we have

\corollary{$\H_G(\C)$ equipped with $\L$ is a conformal vertex algebra.}

The theorem actually holds in slightly greater generality. We do not need to assume that the weights on $\cA$ are nonnegative (this being part of the definition of a half $O(\gs\gg)$ TVA), since the proof never uses the positivity of the weight structure.

The theorem can fail if the center of $G$ has positive dimension. For then $\kappa_G(\L)\circ_n$ will act trivially on $a=\gamma^{\xi_i}\gamma^{\xi_i}\otimes 1$ where the index $i$ here is summed over an orthonormal basis $\xi_i$ of the center of $\gg$. Thus if $a$ represents a nonzero class, then $\kappa_G(\L)$ cannot be a conformal structure in cohomology in this case.

As a corollary, the theorem also gives another proof of Theorem \VirasoroVanishII, since any vertex algebra with a zero conformal structure $L=0$ is necessarily a commutative algebra with only weight zero elements, because $L\circ_1$ must then act by zero. Here is another
immediate consequence.

\corollary{For any $G$ manifold $M$, $G$ semisimple, $\H_G(\cQ'(M))=H_G(M)=\H_G(\cQ(M))$ iff $\kappa_G(\L)=0$.}



\theorem{If $M$ contains a $G$ fixed point, then the chiral Chern-Weil map $\kappa_G:\H_G(\C)\ra\H_G(\cQ'(M))$ is injective. If, moreover, $G$ is semisimple, then the chiral equivariant cohomology is a vertex algebra with a nontrivial conformal structure.}
\proof\thmlab\FixedPointTheorem
By Theorem \FunctorialityQPrime, $\cQ'$ is a contravariant functor from $G$ manifolds to $\gs\gg[t]$ algebras. Since $M$ has a fixed point $pt$, we have equivariant maps $pt\ra M\ra pt$. They induce
$\gs\gg[t]$ algebra homomorphisms $\C\ra\cQ'(M)\ra\C$, whose composition is the identity. So we have the maps
$$
\H_G(\C)\ra\H_G(\cQ'(M))\ra\H_G(\C).
$$
Since the composition is the identity, it follows that the first map, which is $\kappa_G$, is injective. $\Box$

\remark{The proof points to an interesting and much more general fact: any $\gs\gg[t]$ algebra $\cA$ that admits $\C$ as a quotient necessarily has a $G$ chiral equivariant cohomology $\H_G(\cA)$ that contains a copy of $\H_G(\C)$ via the chiral Chern-Weil map. In particular, if $\cF$ is any functor from the category of $G$ spaces to the category of $\gs\gg[t]$ algebras such that $\cF(pt)=\C$, then $\H_G(\cF(M))$ necessarily contains a copy of $\H_G(\C)$ whenever $M$ has a fixed point. The preceding theorem is the special $\cF=\cQ'$.}



\subsec{Chiral equivariant cohomology of $\cQ(M)$}

In this subsection, we give some examples of this cohomology.

Let $f:\cA\ra\cB$ be an $\gs\gg[t]$ algebra homomorphism. Then we have an induced vertex algebra homomorphism on cohomology
$$
f:\H_G(\cA)\ra\H_G(\cB).
$$
In particular, any ``characteristic'' classes that die in $\H_G(\cA)$ will also die in $\H_G(\cB)$. For example, we know the chiral Chern-Weil map $\kappa_G$ is natural, so that the chiral Chern-Weil map into $\H_G(\cB)$ must factor through that into $\H_G(\cA)$. Thus if $G$ is semisimple and if the class $\kappa_G(\L)$ dies in $\H_G(\cA)$, then it must also die in $\H_G(\cB)$. By the same token, any characteristic class that survives in $\H_G(\cB)$ must also survive in $\H_G(\cA)$. Thus the natural factorization of the chiral Chern-Weil map sometimes allows us to detect the vanishing or nonvanishing of certain cohomology classes for one algebra or another.

For any $G$ manifold $M$, the inclusion map $i:\cQ'(M)\hra\cQ(M)$ is clearly an $\gs\gg[t]$ algebra homomorphism. Thus we can apply the above in this geometric setting. Thus if $G$ is semisimple and if $\kappa_G(\L)$ dies in $\H_G(\cQ'(M))$, then it will die in $\H_G(\cQ(M))$ as well. Examples of this include $M=G/T$ where $T$ is a torus subgroup of $G$, as we shall see later.

There are, however, important examples where this factorization cannot detect the vanishing of classes in $\H_G(\cQ(M))$. We now discuss one such example. We will consider linear representations of $G$, viewed as $G$ manifolds. In particular, there is always at least one fixed point, the zero vector in such a representation $V$. Thus by Theorem \FixedPointTheorem, $\H_G(\cQ'(V))$ always contain nontrivial classes with positive weights.

\theorem{Suppose $G$ is semisimple. For any finite dimensional faithful representation $V$ of $G$, the cohomology $\H_G(\cQ(V))$ contains no classes with positive weights.}
\thmlab\LinearRepTheorem

This is a precursor of a theorem about free actions on a manifold, which we shall consider in the next section. We now prepare for the proof. So, we assume that {\it $V$ is a faithful representation of $G$, which is semisimple.} The following lemma is an essential idea that will reappear later. 

Associated to the representation $\rho:\gg\ra End~V$, there is a symmetric invariant bilinear form $\bra\xi,\eta\ket=Tr(\rho(\xi)\rho(\eta))$ on $\gg$. It is nondegenerate since $V$ is faithful. 
For simplicity, we identify $\gg\equiv\gg^*$ via $\bra,\ket$, and we fix an orthonormal basis $\xi_i$ of $\gg$. As usual, we denote by $(\xi,\eta)\mapsto L_\xi+\iota_\eta$ the $O(\gs\gg)$ structure on a given $O(\gs\gg)$ algebra $\cA$. The concrete example we have in mind is  $\cA=\cQ(V)$.

\lemma{Suppose that $(\cA,d_\cA)$ is an $O(\gs\gg)$ algebra with no elements of negative weight and that there is a $\gg$ module homomorphism $\gg\ra\cA[1]$, $\xi\mapsto\Gamma_\xi$, such that for $\xi,\eta,\nu\in\gg$, (1) $\iota_\xi\circ_1\Gamma_\eta=0$, (2) $\iota_\xi\circ_0(\iota_\eta\circ_0\Gamma_\nu)=0$, and (3) $L_\xi\circ_1\Gamma_\eta=\bra\xi,\eta\ket$ hold. Then 
$$
u=\beta^{\xi_i}\partial c^{\xi_i}\otimes1-d\alpha
$$ 
is chiral horizontal in $\cW(\gg)\otimes\cA$, where $\alpha=\beta^{\xi_i}\otimes\Gamma_{\xi_i}-\beta^{\xi_i}c^{\xi_j}\otimes\iota_{\xi_j}\circ_0\Gamma_{\xi_i}$, $d=d_\cW+d_\cA$.}
\proof
We have
$$
(b^\xi+\iota_\xi)\circ_p\alpha=
\beta^{\xi_i}\otimes\iota_\xi\circ_p\Gamma_{\xi_i}-\delta_{p,0}\bra\xi,\xi_j\ket\beta^{\xi_i}\otimes\iota_{\xi_j}
\circ_0\Gamma_{\xi_i}+\beta^{\xi_i}c^{\xi_j}\otimes\iota_\xi\circ_p\iota_{\xi_j}
\circ_0\Gamma_{\xi_i}.
$$
This has negative weight when $p>1$, so it must be zero. It is also zero for $p=1,0$ by assumptions (1) and (2). Thus $\alpha$ is chiral horizontal.

Next, we have
\eqn\dumb{
(b^\xi+\iota_\xi)\circ_pu=\beta^\xi\otimes1~\delta_{p,1}-(\Theta_\cW^\xi+L_\xi)\circ_p\alpha.
}
Here we have used the identity $[d,b^\xi+\iota_\xi]=\Theta_\cW^\xi+L_\xi$ and the fact that the $\alpha$ is chiral horizontal. Since the $\Gamma_\xi$ transform in a (co)adjoint representation of $G$, $\alpha$ is clearly is $G$ invariant. So \dumb~is zero for $p=0$. For $p>1$, \dumb~is also zero because it has negative weight. 
So it remains to show that for $p=1$, the second term on the right side of \dumb~cancels the first term $\beta^\xi\otimes1$. We have
$$
(\Theta_\cW^\xi+L_\xi)\circ_1\alpha=\beta^{\xi_i}\otimes L_\xi\circ_1\Gamma_{\xi_i}-
\beta^{\xi_i}c^{\xi_j}\otimes L_\xi\circ_1\iota_{\xi_j}\circ_0\Gamma_{\xi_i}
$$
where we have used the fact that $\Theta_\cW^\xi\circ_1\beta^{\xi_i}c^{\xi_j}=0$ (no double contraction.)
The first term on the right side becomes $\beta^\xi\otimes1$ by assumption (3). Finally, we have
$$\eqalign{
L_\xi\circ_1\iota_{\xi_j}\circ_0\Gamma_{\xi_i}&=\iota_{\xi_j}\circ_0L_\xi\circ_1
\Gamma_{\xi_i}-(\iota_{\xi_j}\circ_0L_\xi)\circ_1\Gamma_{\xi_j}\cr
&=\iota_{\xi_j}\circ_0\bra\xi,\xi_i\ket-\iota_{[\xi_j,\xi]}\circ_1\Gamma_{\xi_i}=0.
}$$
So, $(\Theta_\cW^\xi+L_\xi)\circ_1\alpha=\beta^\xi\otimes 1$. This completes the proof. $\Box$

\corollary{Under the assumptions of the preceding lemma, $\H_G(\cA)$ has no classes with positive weights.}
\proof
By the preceding lemma, $a=\Theta_\cS^{\xi_i}b^{\xi_i}\otimes1+u$ is also chiral horizontal in $\cW(\gg)\otimes\cA$. But we have
$$
da=d_\cW(\Theta_\cS^{\xi_i}b^{\xi_i}+\beta^{\xi_i}\partial c^{\xi_i})\otimes1\equiv\L\otimes1
$$
which we know is also chiral horizontal. It follows that $a$ must be chiral invariant, implying that $da=\L\otimes 1\equiv\kappa_G(\L)$ is zero in $\H_G(\cA)$. By the corollary to Theorem \ConformalStructureTheorem, $\H_G(\cA)$ has no classes of positive weights. $\Box$

To prove Theorem \LinearRepTheorem, we shall construct a $\gg$ module homomorphism $\gg\ra\cA=\cQ(V)$ with the property prescribed in the preceding lemma.
First we work out explicitly the $O(\gs\gg)$ structure induced on $\cA$ by the $G$ action on the manifold $V$. Fix a set of linear coordinates $\gamma^i\equiv\gamma^{x_i'}$ corresponding to a basis $x_i={\partial\over\partial\gamma^i}$ of $V$. For $\xi\in\gg$, the vector field $X_\xi$ at the point $v=\gamma^i x_i$ is given by
$$
{d\over dt}|_{t=0}e^{t\xi}v=-\rho(\xi)v=-\xi_{ij}\gamma^i x_j
$$
where $\xi_{ij}$ is the matrix of $\rho(\xi)\in End~V$. (Note the sign difference between the $G$ action on space $V$ and its action on functions.) Since $d_\cQ=(\beta^ic^i)(0)$, this shows that the $O(\gs\gg)$ structure on $\cA$ is given by
$$
\iota_\xi=-\xi_{ij}\gamma^ib^j,~~~L_\xi=d_\cQ\iota_\xi=-\xi_{ij}(:\gamma^i\beta^j:+:c^ib^j:)
\equiv-:\beta^{\rho(\xi)x_i}\gamma^{x_i'}:+:b^{\rho(\xi)x_i}c^{x_i'}:.
$$
Note that, formally, $L_\xi$ coincides with $\Theta^\xi_{\cW(V)}$ given in Lemma 2.16 of \LLI.
That is because $L_\xi$ here acts on the space $\cQ(V)$ which happens to contain $\cW(V)$ as a subalgebra.

{\it Proof of Theorem \LinearRepTheorem:}
For each simple summand $\gh$ in $\gg$, we define
$$
\gh\ra\cA=\cQ(V), ~~~\xi\mapsto\Gamma_\xi=:\beta^{\rho(\xi)x_i}\gamma^{x_i'}:.
$$
This defines a $\gg$ module homomorphism $\gg\ra\cA[1]$. By the preceding lemma, it suffices to check conditions (1)-(3) there. Condition (1) holds because the OPE of $\iota_\xi$ and $\Gamma_\eta$ has no double contraction. Condition (2) holds because $\iota_\eta\circ_0\Gamma_\nu$ has the shape $\gamma b$, which commutes with $\iota_\xi$. Finally, for $\eta\in\gh$, a simple summand in $\gg$, we have
$$
L_\xi\circ_1 \Gamma_\eta=-(:\beta^{\rho(\xi)x_i}\gamma^{x_i'}:)\circ_1
(:\beta^{\rho(\eta)x_j}\gamma^{x_j'}:)=\bra\xi,\eta\ket.
$$
Thus condition (3) holds. $\Box$

\remark{The choice of the $\Gamma_\xi$ above is not unique. The choice $:b^{\rho(\xi)x_i}c^{x_i'}:$ would also work.}

\newsec{Vanishing Theorem and Locally Free Actions}

We shall study the chiral equivariant cohomology of manifolds with a locally free $G$ action. We begin by giving the chiral analogue of the classical notion of $G^*$ algebras with property C. The latter is an abstraction of the $G^*$ algebra $\Omega(M)$ when the $G$ action on $M$ is locally free, i.e. every isotropy group is finite.

In this section all $\gs\gg[t]$ algebras $\cA$ are assumed to be contained in some $O(\gs\gg)$ algebra $\cA'$.
As usual, we denote the $O(\gs\gg)$ structure by $(\xi,\eta)\mapsto L_\xi+\iota_\eta$ with $L_\xi=d\iota_\xi$. If $\cA$ is, in addition, a half $O(\gs\gg)$ TVA, then we denote its chiral contracting homotopy by $g$ and conformal structure $L=dg$. When the role of $\cA$ needs emphasis, we will use superscripts or subscripts such as $g^\cA, d_\cA$. {\it We shall further assume that the group $G$ acts on $\cA$ as a group of automorphisms in a way compatible with the half $O(\gs\gg)$ structure.} Namely for $h\in G$, $a\in\cA$, the following relations hold (cf. p. 17 of \GS):
$$
{d\over dt} exp(t\xi)|_{t=0}~a=L_\xi\circ_0 a,~~~
hL_\xi=L_{Ad(h)\xi},~~~
h\iota_\xi=\iota_{Ad(h)\xi},~~~h\circ d\circ h^{-1}=d.
$$
Note that $h\circ\hat L_\xi\circ h^{-1}=\hat L_{Ad(h)\xi}$, and likewise for $\iota_\xi$.
Note that the second relation is a consequence of the third and fourth ones. These relations will be needed in order to average over the group $G$. As in \GS, in order to make sense of differentiation along a curve in $\cA$, we can either assume an appropriate topology on $\cA$ (e.g. $\cQ(M)$), or assume that $\cA$ is $G$ finite (e.g. $\cW(\gg)$), or impose the first condition on $G$ finite elements $a$ only. 

\subsec{Chiral free algebras}

In \LLI, we introduced the notion of a $\cW(\gg)$ algebra, which is a chiral analogue of a $G^*$ algebra with property C \GS. We now introduce another version that better captures the geometric aspect of property C. 

\definition{An $\gs\gg[t]$ algebra $\cA$ is said to be chiral free if there is a linear map $\gg^*\ra\cA^1[0]$, $\xi'\mapsto\theta_{\xi'}$, such that in the $O(\gs\gg)$ algebra $\cA'\supset\cA$ we have
$$
\iota_\xi(z)\theta_{\xi'}(w)\sim\bra\xi',\xi\ket(z-w)^{-1},~~~
L_\xi(z)\theta_{\xi'}(w)\sim\theta_{ad^*(\xi)\xi'}(w)(z-w)^{-1}. 
$$
We shall call $\theta_{\xi'}$ the connection forms of $\cA$. If $\cA$ is, in addition, a half $O(\gs\gg)$ TVA with chiral contracting homotopy $g$, then we put
$$
\Gamma_{\xi'}=g\circ_0\theta_{\xi'}\in\cA^0[1].
$$}


\lemma{The $G$ action on a manifold $M$ is locally free, i.e. the stabilizer of each point is finite, iff the half $O(\gs\gg)$ TVA $\cQ'(M)$ is chiral free.}
\proof
Recall that the $G$ action on $M$ is locally free iff there exists a map $\gg^*\ra\Omega^1(M)$, $\xi'\mapsto\theta_{\xi'}$, such that $\iota_\xi\theta_{\xi'}=\bra\xi,\xi'\ket$ in $\Omega(M)$. This is equivalent to the OPE relation $\iota_\xi(z)\theta_{\xi'}(w)\sim\bra\xi',\xi\ket(z-w)^{-1}$ in $\cQ(M)$. By suitably averaging over $G$, we can assume that in this case the one forms $\theta_{\xi'}$ transform in the coadjoint representation of $G$, so that $L_\xi\theta_{\xi'}=\theta_{ad^*(\xi)\xi'}$ in $\Omega^1(M)$. This relation is equivalent to $L_\xi(z)\theta_{\xi'}(w)\sim\theta_{ad*(\xi)\xi'}(w)(z-w)^{-1}$ in $\cQ(M)$. $\Box$ 

\remark{The preceding lemma holds for the $O(\gs\gg)$ TVA $\cQ(M)$ as well.}

\lemma{Let $\cA$ be an $\gs\gg[t]$ algebra which is $G$ chiral free, and let $\bra\theta\ket$ be the subalgebra generated by the $\theta_{\xi'}$ in $\cA$. Then $\cA$ is free as a $\bra\theta\ket$-module generated by $\cA_{hor}$. In other words, every $a\in\cA$ can be uniquely expressed as $a=\sum p_\alpha a_\alpha$ where $p_\alpha\in\C[\theta_{\xi_i'}(-p-1)|i=1,..,n;~p=0,1,...]$ and the $a_\alpha\in\cA_{hor}$ form a given basis of $\cA_{hor}$. If in addition $\cA$ is abelian, then $\cA\cong\bra\theta\ket\otimes\cA_{hor}$ as vertex algebras.}
\proof\thmlab\Gchiralfree
The second assertion is an immediate consequence of the first.
Since $\iota_\xi(z)\theta_{\xi'}(w)\sim\bra\xi,\xi'\ket(z-w)^{-1}$, we have
$$
[\iota_\xi(p),\theta_{\xi'}(q)]=\bra\xi,\xi'\ket\delta_{p+q,-1}.
$$
For $p\geq0$, $q<0$, these Fourier modes generate a Clifford algebra acting on $\cA$.
Now the first assertion follows from a standard result from the theory of Clifford algebras. $\Box$ 

This lemma is analogous to a result for classical $G^*$ algebras (see Theorem 3.4.1 of \GS.)

\corollary{If $G=T$ is abelian and $\cA$ is $T$ chiral free, then
$\cA_{inv}$ is free as a $\bra\theta\ket$ module generated by $\cA_{bas}$. If, in addition, $\cA$  is abelian, then $\cA_{inv}\cong\bra\theta\ket\otimes\cA_{bas}$ as vertex algebras.}
\proof
Since $T$ is abelian, the $L_\xi$ commute with the $\iota_\xi$ and the $\theta_{\xi'}$. Thus the same Clifford algebra above acts on $\cA_{inv}$. So the same argument applies to this module. $\Box$


For simplicity, we identify $\gg$ with $\gg^*$ using a $G$ invariant form and let $\xi_i$ be an orthonormal basis of $\gg$.

\lemma{Suppose $\cA$ is a half $O(\gs\gg)$ TVA which is chiral free as above. Then $\Gamma_{\xi}=g\circ_0\theta_{\xi}$ are chiral horizontal and we have
$$\eqalign{
&L_\eta\circ_0\Gamma_\xi=\Gamma_{[\eta,\xi]},~~~
L_\eta\circ_1\Gamma_\xi=\bra\eta,\xi\ket,\cr
&\iota_\eta\circ_0\Gamma_\xi=0,~~~
\iota_\eta\circ_1\Gamma_\xi=0,~~~
\iota_\eta\circ_0 d\Gamma_\xi=\Gamma_{[\eta,\xi]},~~~
\iota_\eta\circ_1d\Gamma_\xi=\bra\eta,\xi\ket.
}$$}
\proof\thmlab\GammaTransform
By Lemma \gLGinvariance, $g$ is $G$ invariant. Since the connection forms $\theta_{\xi}$ transform
in the (co)adjoint representation, so do the $\Gamma_\xi$. This proves our first equality.

Since $\iota_\eta\circ_1\theta_\xi=0=\iota_\eta\circ_1 d\theta_\xi$ because $\cA$ is assumed to have no negative weights, and since $g,\iota_\eta$ commute (by definition of a half $O(\gs\gg)$ TVA), we have
$$
L_\eta\circ_1\Gamma_\xi=L_\eta\circ_1(g\circ_0\theta_\xi)=d\iota_\eta\circ_1(g\circ_0\theta_\xi)=-\iota_\eta\circ_1d(g\circ_0\theta_\xi)=-\iota_\eta\circ_1(L\circ_0\theta_\xi).
$$
The last expression is equal to $-\iota_\eta\circ_1\partial\theta_\xi=\iota_\eta\circ_0\theta_\xi=\bra\eta,\xi\ket$, proving our second equality.


Our third and fourth equalities follow from the fact that $g,\iota_\eta$ commute and that $g\circ_0\bra\eta,\xi\ket=0$.
Applying $d$ to the third equality, we get
$L_\eta\circ_0\Gamma_\xi=\iota_\eta\circ_0 d\Gamma_\xi$, proving our fifth equality.
Likewise, applying $d$ to the fourth equality yields our sixth equality. $\Box$


\subsec{Vanishing theorems}

\theorem{Suppose $(\cA,d_\cA)$ is a half $O(\gs\gg)$ TVA which is chiral free as above. Then
$\kappa_G(\gamma^{\xi_i}\partial\gamma^{\xi_i})=0$.}
\proof\thmlab\gammagammatrivial
We first explain the idea behind this. We already know that $$d_\cW(\gamma^{\xi_i}\partial c^{\xi_i}\otimes1)=\gamma^{\xi_i}\partial\gamma^{\xi_i}\otimes1.$$ So to trivialize the right side in $\H_G(\cA)$, we must add an element of $\cW\otimes\cA$ to $\gamma^{\xi_i}\partial c^{\xi_i}\otimes1$, so as to make it $G$ chiral basic. But the modified element $\alpha$ must have the same image under $d$. So whatever we add, it must be $d$ closed. But every $d$ closed element in $\cW\otimes\cA$ is $d$ exact. So we seek $\alpha$ of the form $\alpha=\gamma^{\xi_i}\partial c^{\xi_i}\otimes1+da$. For a half $O(\gs\gg)$ TVA which is chiral free, $a=-\gamma^{\xi_i}\otimes\Gamma_{\xi_i}$ provides just the right modification. We now give the argument. 

We claim that
$$
\alpha=\gamma^{\xi_i}\partial c^{\xi_i}\otimes1-(d_\cW+d_\cA)(\gamma^{\xi_i}\otimes\Gamma_{\xi_i})
$$
is chiral horizontal in $\cW\otimes\cA$, i.e. it commutes with $b^\xi\otimes1+1\otimes\iota_\xi$, and that
$$
(d_\cW+d_\cA)\alpha=\gamma^{\xi_i}\partial\gamma^{\xi_i}\otimes 1.
$$
Since $\alpha$ has weight one, it is enough to see that $(b^\xi\otimes 1+1\otimes\iota_\xi)\circ_p\alpha=0$ for $p=0,1$.
This is a straightforward computation using the preceding lemma and the identity $[d_\cA,\iota_\xi]=L_\xi$. Combining $d_\cW c^\xi=-\half c^{[\xi_j,\xi]}c^{\xi_j}+\gamma^\xi$ and $d_\cW\gamma^\xi=\gamma^{[\xi_j,\xi]}c^{\xi_j}$, we get $d_\cW(\gamma^{\xi_i}\partial c^{\xi_i})=\gamma^{\xi_i}\partial\gamma^{\xi_i}$. 

We have shown that $\alpha$ is chiral horizontal and $(d_\cW+d_\cA)\alpha=\gamma^{\xi_i}\partial\gamma^{\xi_i}\otimes 1$. But the right side of this is also obviously chiral horizontal. So this identity shows that $\alpha$ is chiral invariant, hence chiral basic.  Thus $\gamma^{\xi_i}\partial\gamma^{\xi_i}\otimes 1$ is $d_\cW+d_\cA$ exact in $(\cW\otimes\cA)_{bas}$. $\Box$

\corollary{If the $G$ action on the manifold $M$ is locally free, then $\kappa_G(\gamma^{\xi_i}\partial\gamma^{\xi_i})=0$ in $\H_G^4(\cQ'(M))[1]$.}

In the next section, we shall improve this result and show that if the $G$ action is locally free, then there are no positive weight classes at all. In fact, we will prove this for any half $O(\gs\gg)$ TVA which is chiral free, of which $\cA=\cQ'(M)$ is one example when the $G$ action on $M$ is locally free.
This can be thought of as a chiral analogue of the classical result \GS~that if $A$ is a $G^*$ algebra with property C, then $H_G(A)=H_{bas}(A)$. By Theorem \VirasoroVanish, $\H_{bas}(\cA)=H_{bas}(\cA[0])$ for a half $O(\gs\gg)$ TVA $\cA$. Thus the result says that if $\cA$ is a half $O(\gs\gg)$ TVA which is chiral free, then $\H_G(\cA)=\H_{bas}(\cA)$.

%

We begin with an important case.

\theorem{Suppose $G$ is semisimple. If $(\cA,d)$ is a half $O(\gs\gg)$ TVA which is chiral free, then $\kappa_G(\L)=0$. In particular, $\H_G(\cA)$ has no classes with positive weights.}
\proof\thmlab\SemisimpleChiralFree
The idea is quite similar to that in the preceding theorem. Recall that
$$
\L=d_\cW(\Theta_\cS^{\xi_i}b^{\xi_i}+\beta^{\xi_i}\partial c^{\xi_i})\in\cW(\gg).
$$
To trivialize $\L\otimes 1$ in $\H_G(\cA)$, we must add a $d=d_\cW+d_\cA$ exact element in $\cW\otimes\cA$ to $(\Theta_\cS^{\xi_i}b^{\xi_i}+\beta^{\xi_i}\partial c^{\xi_i})\otimes1$, so as to make the modified element chiral horizontal. Again, we find that adding $-d(\beta^{\xi_i}\otimes\Gamma_{\xi_i})$ provides just the right modification.
Checking that the modified element
$$
(\Theta_\cS^{\xi_i}b^{\xi_i}+\beta^{\xi_i}\partial c^{\xi_i})\otimes1-d(\beta^{\xi_i}\otimes\Gamma_{\xi_i})
$$
is chiral horizontal is a straightforward calculation using Lemma \GammaTransform. $\Box$

As a special case, we have

\theorem{Suppose $G$ is semisimple. If the $G$ action on $M$ is locally free, then $\kappa_G(\L)=0$ in $\H_G(\cQ'(M))$. In particular, $\H_G(\cQ'(M))$ has no classes with positive weights.}
\thmlab\SemisimpleLocallyFree

Later, we shall see that the semisimplicity assumption can be dropped.

\def\q{{\bf q}}
\subsec{The class $\kappa_G(\q)$ and locally free $T$ action}

Fix a half $O(\gs\gg)$ TVA and put
$$
\q=\gamma^{\xi_i}\partial\gamma^{\xi_i}\in\H_G^4(\C)[1].
$$
Under the assumption that $\kappa_G(\q)=0$ in $\H_G(\cA)$, can we reconstruct the connection forms $\theta_{\xi'}$, and hence show that $\cA$ is chiral free? We shall see that in the geometric setting $\cA=\cQ'(M)$ the answer is affirmative when $G$ is {\it abelian}, but not so in general.

A direct attack on this problem in the chiral Weil model seems to require quite substantial information about the linear $G$ action on the weight one subspace of $\cW\otimes\cQ'(M)$. This strategy amounts to showing that if there exists a point $o\in M$ whose stabilizer group $G_0$ has positive dimension, then there is no way to construct connection forms $\theta_\xi$, and from it, the chiral horizontal element $\alpha$ above, to trivialize the class $\kappa_G(\q)$. The difficulty here is that while the hypothetical $\alpha$ is required to be $G$ invariant, the chiral horizontal condition involves operators which are not $G$ invariant, but rather transform in the adjoint representation of $G$. 

We now discuss an alternative approach. We will first work in the general algebraic setting of a half $O(\gs\gg)$ TVA $\cA$ contained in an $O(\gs\gg)$ TVA $\cA'$, but later specialize to the geometric setting $\cA=\cQ'(M)$ of a $G$ manifold.
Our approach is to use a kind of Cartan model to circumvent the chiral horizontal condition, while paying the price of complicating the differential. However, this gives rise to another problem: $\cW\otimes\cA$ naively has no Cartan model because $\cA$ is {\it not} an $O(\gs\gg)$ algebra; it is merely an $\gs\gg[t]$ algebra. In particular, the chiral Mathai-Quillen map $\Phi=e^{\phi(0)}$ is not well-defined on $\cW\otimes\cA$ because the vertex operator
$$
\phi=c^{\xi_i}\otimes\iota_{\xi_i}
$$
does not belong in $\cW\otimes\cA$. One might then try to modify $e^{\phi(0)}$ by, say, truncating the zeroth mode $\phi(0)$ to $\phi_+=\sum_{n\geq0}c^{\xi_i}(-n-1)\otimes\iota_{\xi_i}(n)$ so that it now operates on $\cW\otimes\cA$. But the problem is that this truncated operator $\phi_+$ is {\it not} a derivation of the circle products. (E.g. $\phi_+[(b\otimes c)\circ_{-1}(\partial b\otimes \partial c)]\neq\phi_+(b\otimes c)\circ_{-1}(\partial b\otimes \partial c)+(b\otimes c)\circ_{-1}\phi_+(\partial b\otimes \partial c)$.) But here is a crucial observation:
$\cW\otimes\cA$ is a differential subalgebra of an algebra $\cW\otimes\cA'$, which has an $O(\gs\gg)$ structure and on which the chiral Mathai-Quillen map $\Phi$ is well-defined. In particular,
the subcomplex
$$
(\cW\otimes\cA)_{bas}\subset(\cW\otimes\cA')_{bas}
$$
is mapped to some subcomplex of $\Phi(\cW\otimes\cA')_{bas}=(\cW_{hor}\otimes\cA')_{inv}$, which we can analyze.

\lemma{Let $\cB$ be the vertex subalgebra of $\cA'$ generated by $\cA$ and the vertex operators $\iota_\xi$. Then $\Phi(\cW\otimes\cA)_{bas}\subset(\cW_{hor}\otimes\cB)_{inv}$.}
\proof
Note that the algebra $\cW\otimes\cB$ need not be a complex, but it is a module over $\gs\gg[t]$.
It is a subspace of $\cW\otimes\cA'$ where the action of $O(\gs\gg)$ is defined;
$(\cW_{hor}\otimes\cB)_{inv}$ is simply the intersection $(\cW_{hor}\otimes\cA')_{inv}\cap(\cW\otimes\cB)$.

Clearly the operator $\Phi$ operates on the vertex algebra $\cW\otimes\cB$ containing $\cW\otimes\cA$. In particular, we have
$$
\Phi(\cW\otimes\cA)_{bas}\subset(\cW_{hor}\otimes\cA')_{inv}\cap(\cW\otimes\cB).~~~~\Box
$$

Note that
$$
\Phi(\gamma^{\xi_i}\partial\gamma^{\xi_i}\otimes1)=\gamma^{\xi_i}\partial\gamma^{\xi_i}\otimes 1
$$
and that the differential on the complex $\Phi(\cW\otimes\cA)_{bas}$, restricted to weight one, is
$$\eqalign{
d_G=\Phi(d_\cW+d_\cA)\Phi^{-1}&=d_\cW+d_\cA-\gamma^{\xi_i}(0)\otimes\iota_{\xi_i}(-1)-\gamma^{\xi_i}(-1)\otimes\iota_{\xi_i}(0)-\gamma^{\xi_i}(-2)\otimes\iota_{\xi_i}(1)\cr
&+c^{\xi_i}(0)\otimes L_{\xi_i}(-1)+c^{\xi_i}(-1)\otimes L_{\xi_i}(0)+c^{\xi_i}(-2)\otimes L_{\xi_i}(1).
}$$

We now study the shape of an element $a$ such that
\eqn\dumb{
d_Ga=\q\otimes1.
}
Since $\cW=\cW(\gg)$ is generated by weight and degree homogeneous vertex operators $b,c,\beta,\gamma$,
we can decompose each $\cW^p[n]$ into a direct sum of finite dimensional subspaces labeled by the shapes of the monomials that span them. For example, if $V$ is a given subspace of $\cB$, we call the subspace of $\cW^4[1]\otimes V$ spanned by all monomials $\gamma^\xi\partial\gamma^\eta\otimes v$ the {\it shape} $\gamma\partial\gamma\otimes V$. Assuming \dumb, we now examine how each operator appearing in $d_G$ and each component of $a$ in a given shape can contribute to $\q\otimes 1\in\gamma\partial\gamma\otimes\cB^0[0]\subset\cW^4[1]\otimes\cB^0[0]$. 

First we consider all possible shapes occuring in $a\in (\cW_{hor}\otimes\cB)[1]$ which can contribute to $\q\otimes1$ in \dumb. By consideration of degrees on both $\cW$ and $\cB$, and degrees of each term in $d_G$, the possible shapes in $\cW_{hor}[1]\otimes\cB[0]$ which can contribute to $\q\otimes1$ in \dumb~are
$$
\beta\gamma\gamma\otimes\cB^1[0],
\beta\gamma\otimes\cB^3[0],\beta\otimes\cB^5[0],
\underline{\partial\gamma\otimes\cB^1[0]},\underline{b\gamma\gamma\otimes\cB^0[0]},
b\gamma\otimes\cB^2[0], b\otimes\cB^4[0].
$$
The possible shapes in $\cW_{hor}[0]\otimes\cB[1]$ which can contribute are
$$
\gamma\gamma\otimes\cB^{-1}[1],\underline{\gamma\otimes\cB^1[1]},
1\otimes\cB^3[1].
$$
Let $\pi:\cW\otimes\cB\ra\gamma\partial\gamma\otimes\cB^0[0]$ be the projection which sends all $\cW^p[m]\otimes\cB^q[n]$ to zero unless $p=4,q=0,m=1,n=0$, in which case $\pi$ sends all shapes $\lambda\otimes\cB^0[0]$ to zero unless $\lambda=\gamma\partial\gamma$. In the last case, $\pi$ is the identity map. Note that this projection map is a $G$-module homomorphism. Then by weight and degree consideration of those shapes listed above, we see that only those underlined survive in $\pi d_G a$. Thus \dumb~becomes
$$
\q\otimes1=\pi d_G a=\pi d_G(a[b\gamma\gamma]+\partial\gamma^{\xi_i}\otimes f_{\xi_i}+\gamma^{\xi_i}\otimes g_{\xi_i})
$$
where $a[b\gamma\gamma]$ is the component of $a$ of the shape $b\gamma\gamma\otimes\cB^0[0]$, and  $\gg\ra\cB^1[0]=\cA^1[0]$, $\xi\mapsto f_\xi$, and $\gg\ra\cB^1[1]$, $\xi\mapsto g_\xi$, are some $G$ module homomorphisms. Let 
$$
\pi_0:(\gamma\partial\gamma\otimes\cA^0[0])^G\ra
(\gamma\partial\gamma)^G\otimes\cA^0[0]^G
$$ 
be the projection map that sends all subspaces of the form $\gamma\partial\gamma[V]\otimes *$ to zero, except when $V=\C$, in which case $\pi_0$ is the identity map. Here $\gamma\partial\gamma[V]$ is the $G$ isotypic subspace of $\gamma\partial\gamma$ of the irreducible $G$ type $V$.
Since $\q\in(\gamma\partial\gamma)^G$, it follows that $\q\otimes1=\pi_0\pi d_G a$. Since $d_\cW$ is a $G$ homomorphism and $\pi_0\pi d_G b\gamma\gamma=\pi_0\pi d_\cW b\gamma\gamma$, the contribution of this term to $\q\otimes1$ will only come from $(b\gamma\gamma)^G$. By Lemma 7.8 of \LLI, $d_\cW(b\gamma\gamma)^G=0$. Thus we get

\lemma{A necessary condition for $\q\otimes1$ to be $d_G$ exact is that there exist $G$ module homomorphisms $\gg\ra\cB^1[0]$, $\xi\mapsto f_\xi$, and $\gg\ra\cB^1[1]$, $\xi\mapsto g_\xi$, such that
\eqn\dumb{
\q\otimes1=\pi_0(\gamma^{\xi_i}\partial\gamma^{\xi_j}\otimes h_{ij})
}
where $-h_{ij}=\iota_{\xi_i}\circ_0 f_{\xi_j}+\iota_{\xi_j}\circ_1 g_{\xi_i}$.}
\proof
The $f$ term comes from $-\gamma^{\xi_i}(-1)\otimes\iota_{\xi_i}(0)$, while
the $g$ term comes from $-\gamma^{\xi_i}(-2)\otimes\iota_{\xi_i}(1)$ of $d_G$ acting on
$\partial\gamma^{\xi_i}\otimes f_{\xi_i}$, $\gamma^{\xi_i}\otimes g_{\xi_i}$, respectively. 
By weight, degree and shape considerations, after projecting under $\pi_0\pi$, no other terms will contribute to $\q\otimes1$. $\Box$

We now specialize above to the geometric setting where $\cA=\cQ'(M)$, $\cA'=\cQ(M)$, for a $G$ manifold $M$. Thus $\cB$ is the vertex subalgebra generated by $\cA$ and the $\iota_\xi$ in $\cQ(M)$.

\lemma{Suppose there exists a point $o\in M$ whose stabilizer group $G_0$ is infinite. Then there exists a nonzero $\xi\in\gg$ such that, for all $x\in\cB$,
$$
\iota_\xi\circ_p x|_o=0,~~~p\geq0.
$$}
\proof
Since $G$ is compact, $G_0$ has positive dimension. Hence there exists a nonzero element $\xi\in\gg$ such that the vector field $x_\xi$ generated by $\xi$ vanishes at $o$. In particular $\iota_\xi\circ_0 f=\iota_\xi f$ vanishes at $o$. In local coordinates, we have $\iota_\xi=f^\xi_i b^i$ where the $f^\xi_i$ are functions which vanish at $o$. Since $\cB$ is the vertex algebra generated by $\cQ'(M)$ and the vertex operators $\iota_\eta$, it suffices to check our assertion for generators $x$, by Lemma \Associativity. For $x\in\cQ'(M)$, this is clear. For $x=\iota_\eta$, we have $\iota_\xi\circ_p\iota_\eta=0$ for $p\geq0$. $\Box$

A word of caution is in order here. Evaluation of a vertex operator at a point $o\in M$ must be done with care. For a general $x\in\cQ(M)$, $x|_o$ has no meaning. The reason is that in local coordinates $\gamma^i$, even though each $\gamma^i$ is a function vanishing at the origin, we have $\beta^j\circ_0\gamma^i=\delta_{ij}$. When the functions appearing in $x$ are vanishing with respect to a chart, it need not be so with respect to a different chart. Thus vanishing at a point is not a coordinate invariant notion! This phenomenon is closely related to the fact that the sheaf $\cQ$ is not a module over $C^\infty$. However if $x\in\cB$, $x|_o$ makes sense because the restriction of $x$ to any local chart does not involve the vertex operators $\beta^i$. In any case, for what we need, it is enough for the preceding lemma to hold for $x\in\Omega^1(M)$ and $x\in\cB^1[1]$, which can be checked explicitly in any local coordinates.

\theorem{Suppose $G$ is abelian. Then $\kappa_G(\q)=0$ in $\H_G(\cQ'(M))$ iff the $G$ action on $M$ is locally free.}
\proof\thmlab\AbelianFreeAction
We have already proved the \lq\lq if" part.
In this case, $G$ acts trivially on $\cW$ and hence each $\gamma^\xi\partial\gamma^\eta$ is $G$ invariant, and $\pi_0$ is the identity map. Then the preceding lemma yields that $h_{ij}=\delta_{ij}$.
Suppose the $G$ on $M$ is not locally free. Then there is a point $o\in M$ such that the stabilizer group has positive dimension. So there is a unit vector $\xi$ in the Lie algebra of the stabilizer of $o$, such that the vector field $X_\xi$ vanishes at $o$. Then for $\iota_\xi=f^\xi_i b^i$, we have $f^\xi_i|_o=0$. In particular, $\iota_\xi\circ_p x|_o=0$ for all $x\in\cB$ and $p\geq0$, by the preceding lemma. Choose $\xi_i$ so that $\xi_1=\xi$. Then $h_{11}|_o=0$, contradicting $h_{11}=\delta_{11}=1$. $\Box$

We shall see in the next section that the vanishing of $\kappa_G(\q)$ is not enough to guarantee locally free action when $G$ is semisimple.

\newsec{Chiral Equivariant Cohomology of a Quotient Space}

\subsec{$T$-reduction}

Throughout this subsection, let $G,T$ be compact connected Lie groups; $T$ is assumed abelian.

\theorem{If $\cA$ is a half $O(\gs\gg\oplus\gs\gt)$ TVA which is $T$ chiral free, then
the chain map $(\cW(\gt)\otimes\cW(\gg)\otimes\cA)_{G\times T-bas}\hla(\cW(\gg)\otimes\cA_{T-bas})_{G-bas}$ induces an isomorphism 
$$\H_{G\times T}(\cA)\cong\H_G(\cA_{T-bas}).$$}
\thmlab\TReduction

Before the proof, we give a few important consequences.

\corollary{If $M$ is $G\times T$ manifold on which $T$ acts freely, then
$$\H_{G\times T}(\cQ'(M))\cong\H_G(\cQ'(M/T)).$$}
\thmlab\TReductioncor
\proof
Since $T$ acts freely, $\cQ'(M)$ is $T$ chiral free. By Theorem \SheafGReduction, $\cQ'(M)_{T-bas}\cong\cQ'(M/T)$. Now apply the preceding theorem to $\cA=\cQ'(M)$. $\Box$

\theorem{If $\cA$ is a half $O(\gs\gg)$ TVA which is $G$ chiral free, then $\H_G(\cA)$ has no classes with positive weights.}
\proof\thmlab\GeneralChiralFree
Since $G$ is assumed connected, it is of the form $(H\times T)/F$ where $H$ is semisimple, $T$ is a torus, and $F$ is a finite subgroup of $H\times T$. Thus given a half $O(\gs\gg)$ TVA $\cA$,
we can regard it as a half $O(\gs\gh\oplus\gs\gt)=O(\gs\gg)$ TVA on which the subgroup $F\subset H\times T$ acts trivially. So we have
$$
\H_G(\cA)=\H_{H\times T}(\cA).
$$
By the preceding theorem, it suffices to show that $\cA_{T-bas}$ is a half $O(\gs\gh)$ TVA which is $H$ chiral free. For then Theorem \SemisimpleChiralFree, applied to this algebra, yields our assertion. 

Since the actions of $H$ and $T$ commute, the $\gs\gh[t]$ structure on $\cA$ restricts to $\cA_{T-bas}$, i.e. for $\xi\in\gh$, $p\geq0$, the operators $L_\xi\circ_p,\iota_\xi\circ_p$ preserve $\cA_{T-bas}$. By definition, the chiral contracting homotopy $g$ is $G$ chiral horizontal. In particular, for $\eta\in\gt$, $p\geq0$, we have
$$
L_\eta\circ_p g=-\iota_\eta\circ_p dg=-\iota_\eta\circ_p L=\iota_\eta~\delta_{p,1}
$$
where the last equality holds because $\iota_\eta$ is primary of weight one. This shows that $\cA_{T-bas}$ is also preserved by the action of $g\circ_n$, $n\geq0$. Thus the half $O(\gs\gh)$ TVA structure on $\cA$ also restricts to the subalgebra $\cA_{T-bas}$. Finally, we need to show that this structure is $H$ chiral free. By assumption, $\cA$ is $G$ chiral free. So we have $G$ connection forms $\theta_{\xi'}$ parameterized by $\xi'\in\gg=\gh\oplus\gt$. By Lemma \GammaTransform,
we see that the connection forms parameterized by $\xi'\in\gh$ are all $T$ chiral basic because $\gh^*\perp\gt$. In other words, we have $\gh^*\ra\cA_{T-bas}$, $\xi'\mapsto\theta_{\xi'}$, defining $H$ connection forms on $\cA_{T-bas}$. This completes the proof. $\Box$

\corollary{If $M$ is a $G$ manifold on which the $G$ action is locally free, then
$\H_{G}(\cA)$ has no classes with positive weights for $\cA=\cQ'(M)$ or $\cQ(M)$.}

This generalizes Theorem \SemisimpleLocallyFree.
We now proceed to the proof of Theorem \TReduction.

\lemma{Put $\delta=d_{\cW(\gt)}+d_{\cW(\gg)}+d_\cA$, and let $\bra\gamma,c\ket_T$ be the subalgebra of $\cW(\gt)$ generated by the $\gamma^{\eta'},c^{\eta'}$ with $\eta'\in\gt^*$. Then
the chain map $(\cW(\gt)\otimes\cW(\gg)\otimes\cA)_{G\times T-bas}\hla[\bra\gamma,c\ket_T\otimes(\cW(\gg)\otimes\cA_{T-inv})_{G-bas}]_{T-hor}$ induces an isomorphism
$$
\H_{G\times T}(\cA)\cong H\left([\bra\gamma,c\ket_T\otimes(\cW(\gg)\otimes\cA_{T-inv})_{G-bas}]_{T-hor},\delta\right).
$$}
\proof\thmlab\TwoStageLemma
Since $\cW(\gg\oplus\gt)=\cW(\gt)\otimes\cW(\gg)$, the chiral Weil model gives
$$
\H_{G\times T}(\cA)=H\left([\cW(\gt)\otimes(\cW(\gg)\otimes\cA)_{G-bas}]_{T-bas},\delta\right).
$$
We can view $(\cW(\gg)\otimes\cA)_{G-bas}$ as an $\gs\gt[t]$ subalgebra of the $O(\gs\gt)$ algebra $(\cW(\gg)\otimes\cA')_{G-bas}$ with differential $d_{\cW(\gg)}+d_\cA$, and consider the small chiral Weil model of this subalgebra. Here $\cA'$ is some $O(\gs\gg\oplus\gs\gt)$ TVA containing $\cA$. Since $O(\gs\gt)$ acts trivially on $\cW(\gg)$, we have 
$$
[(\cW(\gg)\otimes\cA)_{G-bas}]_{T-inv}=(\cW(\gg)\otimes\cA_{T-inv})_{G-bas}.
$$
Now Theorem \SmallWeilModel~(small chiral Weil model) yields the desired isomorphism. $\Box$

We would like to take further advantage of the $T$ chiral freeness of $\cA$ by using Lemma \Gchiralfree~ to describe the algebra
$$
[\bra\gamma,c\ket_T\otimes(\cW(\gg)\otimes\cA_{T-inv})_{G-bas}]_{T-hor}=
[\bra\gamma,c\ket_T\otimes(\cW(\gg)\otimes\cA)_{G-bas}]_{T-bas}.
$$
Unfortunately, the $T$ connection forms $\theta_{\xi'}$ are almost never $G$ chiral basic because, in general, they don't commute with $\iota_\xi$, $\xi\in\gg$, i.e. they fail to be $G$ chiral horizontal. To circumvent this, we use the $G$ chiral Mathai-Quillen isomorphism $\Phi_G$ to isolate the $G$ chiral horizontal condition to the factor $\cW(\gg)$. Recall that
$$
\Phi_G=e^{(c^{\xi_i'}\otimes\iota_{\xi_i})(0)}
$$
which obviously commutes with the $T$ chiral basic condition. Thus
$$\eqalign{
&\Phi_G[\bra\gamma,c\ket_T\otimes(\cW(\gg)\otimes\cA)_{G-bas}]_{T-bas}
=(\bra\gamma,c\ket_T\otimes\cB)_{T-bas},\cr
&\cB{\br def\over=}\Phi_G(\cW(\gg)\otimes\cA)_{G-bas}\subset(\cW(\gg)_{G-hor}\otimes\cA')_{G-inv}.
}$$
Clearly $\cB$ is an $\gs\gt[t]$ algebra. Since the $T$ connection forms $\theta_{\eta'}\in\cA^1[0]$ are chosen to be invariant under $G$, and since they have weight zero, they are also $G$ chiral invariant. Hence they must lie in $(\cW(\gg)_{G-hor}\otimes\cA')_{G-inv}$. This implies that $\Phi_G^{-1}\theta_{\eta'}=\theta_{\eta'}-c^{\xi_i'}\otimes(\iota_{\xi_i}\circ_0\theta_{\eta'})\in\cW(\gg)\otimes\cA$ is $G$ chiral basic. This shows that $\Phi_G(\Phi_G^{-1}\theta_{\eta'})=\theta_{\eta'}\in\cB$ (but note that $\Phi_G\theta_{\eta'}\neq\theta_{\eta'}$). It follows that $\cB$ is an $\gs\gt[t]$ algebra which is $T$ chiral free. Note that $\bra\gamma,c\ket_T$ is also an $\gs\gt[t]$ algebra which is $T$ chiral free with connection forms $\eta'\mapsto c^{\eta'}$ (and trivial $\gt[t]$ action). 
The differential on $(\bra\gamma,c\ket_T\otimes\cB)_{T-bas}$ is
$$
\Phi_G\delta\Phi_G^{-1}=d_{\cW(\gt)}+d_G,~~~d_G=d_\cW(\gg)+d_\cA-(\gamma^{\xi_i'}\otimes\iota_{\xi_i})(0)+(c^{\xi_i'}\otimes L_{\xi_i})(0)
$$
where $d_G$ acts only on the $\cB$ factor.


We now make a few observations about the complex $\left((\bra\gamma,c\ket_T\otimes\cB)_{T-bas},\Phi_G\delta\Phi_G^{-1}\right)$.

\bu For $\eta\in\gt$, $b^\eta+\iota_\eta$ commutes with $1\otimes\cB_{T-bas}$, because $b^\eta,\iota_\eta$ do separately.

\bu Since $T$ acts trivially on $\bra\gamma,c\ket_T$, $(\bra\gamma,c\ket_T\otimes\cB)_{T-bas}$ is the subalgebra of $(\bra\gamma,c\ket_T\otimes\cB)_{T-inv}=\bra\gamma,c\ket_T\otimes\cB_{T-inv}$ commuting with the $b^\eta+\iota_\eta$.

\bu $\bra\gamma,c\ket_T\otimes\cB_{T-inv}$ is a free module over the algebra $\bra\gamma,c,\theta\ket_T$ on the subspace $1\otimes\cB_{T-bas}$.

\bu For $\eta'\in\gt'$, put
$$
u^{\eta'}=c^{\eta'}+\theta_{\eta'},~~~~v^{\eta'}=c^{\eta'}-\theta_{\eta'}.
$$
Then the $\gamma^{\eta'},v^{\eta'}$ commute with the $b^\eta+\iota_\eta$ and the $u^{\eta'}$; the nonnegative Fourier modes of the $b^\eta+\iota_\eta$ and the negative modes of the $u^{\eta'}$ generate a Clifford algebra.

\bu Thus, $(\bra\gamma,c\ket_T\otimes\cB)_{T-bas}$ is free over $\bra\gamma,v\ket$ on $1\otimes\cB_{T-bas}$.

\lemma{Let $(W^*,\delta_1)$ be an abstract Koszul algebra, i.e. $W^j=\oplus_{j=p+2q}\wedge^p V\bigoplus S^qV$ and $\delta_1$ is the derivation such that $\delta_1(v\otimes 1)=1\otimes v$ for $v\in V$. Let $C$ be a free $W$-module generated by a subspace $B\subset C$. Assume $C$ is equipped with two anticommuting square-zero differentials $\delta_1,\delta_2$ and an inclusion $(W,\delta_1)\hra(C,\delta_1)$ as differential $W$-modules, such that
\eqn\dumb{
\delta_1(wb)=(\delta_1w)b,~~w\in W,~b\in B,~~~\delta_2:W^j B\ra W^jB\oplus W^{j-1}B.
}
(In particular $\delta_1|B=0$.)
Then the inclusion $(B,\delta_2)\hra(C,\delta_1+\delta_2)$ is a quasi-isomorphism.}
\proof
Since $C$ is $W$ free on $B$, we have $C=\bigoplus_j C^j$ where $C^j=W^jB$, and
$$
H^j(C^*,\delta_1)=B\delta_{j,0}.
$$
Now our assertion follows from a standard result in homological algebra \HS. $\Box$

{\it Proof of Theorem \TReduction:} 
Put
$$
C=(\bra\gamma,c\ket_T\otimes\cB)_{T-bas},~~~B=1\otimes\cB_{T-bas},~~~W=\bra\gamma,v\ket.
$$
View $W$ as the Koszul algebra generated by the negative Fourier modes of the $\gamma^{\eta'},v^{\eta'}$, for $\eta\in\gt$, $\eta'\in\gt^*$, and view $C$ as an abstract free module over $W$ on $B$. The differentials are
$$
\delta_1=d_{\cW(\gt)}=(\gamma^{\eta_i'}b^{\eta_i})(0),~~~
\delta_2=d_G=d_\cW(\gg)+d_\cA-(\gamma^{\xi_i'}\otimes\iota_{\xi_i})(0)+(c^{\xi_i'}\otimes L_{\xi_i})(0).
$$
They anticommute, $\delta_1$ kills $B$, and $\delta_1 v^{\eta'}=\gamma^{\eta'}$. Since $(\delta_1+\delta_2)v^{\eta'}$ lies in $C$, it must commute with the $b^\eta+\iota_\eta$ for $\eta\in\gt$. This implies that $\delta_2v^{\eta'}=-\delta_2\theta^{\eta'}$ must commute with $\iota_\eta$. (This can also be verified by direct computation.) Since $\theta^{\eta'}\in\cB$, so is $\delta_2\theta^{\eta'}\in\cB$, because $(\cB,\delta_2)$ is a complex. This shows that $\delta_2v^{\eta'}\in\cB$. Since $\delta_2$ is $T$ invariant and $\delta_2\theta^{\eta'}$ has weight zero, it is also $T$ chiral invariant. But a $T$ chiral invariant element in $\cB$ lies in $B=1\otimes\cB_{T-bas}\equiv\cB_{T-bas}$ iff it commutes with the $\iota_\eta$. This shows that 
$$
\delta_2 v^{\eta'}\in B.
$$
Since $\partial$ commutes with $\delta_2$, we also have $\delta_2\partial^p v^{\eta'}\in B$. Since $\delta_2 v^{\eta'}\in C[0]$ and $C[0]$ is a commutative algebra (with product $\circ_{-1}$), this element commutes with all $\gamma,v$, hence with $W$. It follows that all $\partial$ derivatives of $\delta_2 v^{\eta'}$ also commute with $W$. Clearly $\delta_2\gamma^{\eta'}=0$. So for $w\in W^j$, we have $[\delta_2, w]\in W^{j-1}B$.
Hence we conclude that
$$
\delta_2:W^jB\ra W^jB\oplus W^{j-1}B.
$$
Thus we have verified \dumb~of the preceding lemma. So the map $B\hra C$ induces
$$
H(B,\delta_2)\cong H(C,\delta_1+\delta_2).
$$
We also have
$$
H(B,\delta_2)=H(\cB_{T-bas},d_G)
=H(\Phi_G(\cW(\gg)\otimes\cA_{T-bas})_{G-bas},\Phi_G\delta\Phi_G^{-1})
$$
where the right side is nothing but the chiral Cartan model of $\H_G(\cA_{T-bas})$. To summarize, we have the chain maps
$$
(\cW(t)\otimes\cW(\gg)\otimes\cA)_{G\times T-bas}\hla[\bra\gamma,c\ket_T\otimes(\cW(\gg)\otimes\cA)_{G-bas}]_{T-bas}
{\br\Phi_G^{-1}\over\la}C\hla B{\br\Phi_G\over\la}A
$$
where $A=(\cW(\gg)\otimes\cA_{T-bas})_{G-bas}$. We have shown that each of these chain maps
induces an isomorphism on cohomology. It is clear that their composition is the canonical inclusion $(\cW(t)\otimes\cW(\gg)\otimes\cA)_{G\times T-bas}\hla(\cW(\gg)\otimes\cA_{T-bas})_{G-bas}$. This completes the proof of Theorem \TReduction. $\Box$


\remark{The results proved in this subsection require only that $\cA$ possesses a $T$ chiral free $(\gs\gg\oplus\gs\gt)[t]$ algebra structure, and not a half $O(\gs\gg\oplus\gs\gt)$ TVA structure, because the chiral contracting homotopy $g$ and the conformal structure $L=dg$ are never used.}

\subsec{Chiral connection forms}

Throughout this subsection, $G$ and $H$ will be compact connected Lie groups.

Let $\cA$ be a half $O(\gs\gg\oplus\gs\gh)$ TVA which is $H$ chiral free
and let $\gh^*\ra\theta_{\eta'}\in\cA[0]$ be $H$ connection forms. By averaging over $G$ (since $G$ and $H$ commute), we can assume that the $\theta_{\eta'}$ are $G$ invariant. Since $\theta_{\eta'}$ has weight zero, it is then $G$ chiral invariant.  
{\it We shall further assume that for $\eta'\in\gh^*$}
\eqn\dumbStar{
\Gamma_{\eta'}=g\circ_0\theta_{\eta'}\in\cA_{G-hor},~~~\theta_{\eta'}\in\cA_{G-inv},
}
where $g$ is the chiral contracting homotopy for the half $O(\gs\gg\oplus\gs\gh)$ TVA $\cA$. By definition, $g$ commutes with $\iota_\xi$ for all $\xi\in\gg\oplus\gh$. Thus for $p\geq0$ we have $\iota_\xi\circ_p\Gamma_{\eta'}=g\circ_0(\iota_\xi\circ_0\theta_{\eta'})\delta_{p,0}$,
which a priori need not be zero when $p=0$, unless $\xi\in\gh$. Thus the first condition in \dumbStar~is an additional restriction on $\cA$. We shall call $\Gamma_{\eta'}$ a chiral connection form of $\cA$.

{\it Geometrical example.} In any case, \dumbStar~holds in the geometrical setting where $\cA=\cQ'(M)\subset\cQ(M)$ and $M$ is a $G\times H$ manifold on which the $H$ action is locally free.
In this case, since $g=b^i\partial\gamma^i$ (in local coordinates) and $\theta_{\eta'}\in\Omega^1(M)$ is chosen to be $G$ invariant, it follows that each of the $\Gamma_{\xi'}$ is a sum of $\partial\gamma^i$ with coefficients given by smooth functions. Hence it commutes with $\iota_X$ for any vector field $X$. In particular $\Gamma_{\eta'}$ is $G\times H$ chiral horizontal, and so assumption \dumbStar~holds in this case.

Throughout this subsection,
{\it $\cA$ will be a half $O(\gs\gg\oplus\gs\gh)$ TVA which is chiral free with respect to $G$ and $H$ separately, such that assumption \dumbStar~holds.}

\lemma{The map $\gh^*\ra\cA$, $\eta'\mapsto\Gamma_{\eta'},$ is a homomorphism of $G\times H$ modules. Moreover, for $\xi\in\gg\oplus\gh,\eta'\in\gh^*$,
\eqn\dumb{
L_\xi\circ_1\Gamma_{\eta'}=\iota_\xi\circ_0\theta_{\eta'}\in\cA[0].
}
In particular, for $\xi\in\gt$, $\iota_\xi\circ_0\theta_{\eta'}$ is equal to the scalar $\bra\xi,\eta'\ket$. All three statements hold if the roles of $G$ and $H$ are interchanged.}
\proof\thmlab\GammaGTTransform
The map $\gh^*\ra\cA$, $\eta'\mapsto\theta_{\eta'}$, is clearly an $H$ module homomorphism. Since we have chosen the $\theta_{\eta'}$ to be $G$ invariant (cf. \dumbStar), this is a $G\times H$ module homomorphism. Composing it with $g\circ_0$, we get $\eta'\mapsto\Gamma_{\eta'}$. Since $g$ is $G\times H$ invariant by Lemma \gLGinvariance, our first assertion follows.

For the second statement, we cannot just apply Lemma \GammaTransform~naively here because the connection forms $\theta_{\eta'}$ are with respect to the $H$ action only. 
But the proof of the second equality there carries over to prove \dumb~with only trivial changes. Finally, since the $\theta_{\eta'}$ are $H$ connection forms, the right side of \dumb~ is $\bra\xi,\eta'\ket$, when $\xi\in\gh$. $\Box$


\subsec{$\H_{G\times H}(\cA)$}

Throughout this subsection, $G$ and $H$ will be compact connected Lie groups; $G$ is assumed semisimple. 

We can regard $\cW(\gg)$ as an $O(\gs\gg\oplus\gs\gh)$ algebra in which $O(\gs\gh)$ acts trivially.
In particular, it is automatically $H$ chiral basic. Thus its $G\times H$ chiral basic cohomology is nothing but
$$
\H_{G\times H-bas}(\cW(\gg))=\H_G(\C).
$$
Consider the map
$$
\cW(\gg)\hra\cW(\gg)\otimes\cW(\gh),~~~x\mapsto x\otimes 1,
$$
which is obviously a homomorphism of $O(\gs\gg\oplus\gs\gh)$ algebras. This induces a map on the $G\times H$ chiral basic cohomology. So we have a canonical map
$$
\H_G(\C)\ra\H_{G\times H}(\C),~~~[x]\mapsto[x\otimes 1].
$$
Thus if $\cA$ is any $(\gs\gg\oplus\gs\gh)[t]$ algebra, we can compose this map with the $G\times H$ chiral Chern-Weil map of $\cA$ and get
$$
\H_G(\C)\ra\H_{G\times H}(\cA).
$$
In particular, we can view $\H_{G\times H}(\cA)$ as a module over the vertex algebra $\H_G(\C)$.
Likewise we can do the same for $\H_H(\C)$.

\lemma{For any $(\gs\gg\oplus\gs\gh)[t]$ algebra, the natural map
$$
\H_G(\C)\otimes\H_H(\C)\ra\H_{G\times H}(\cA),~~~[x]\otimes[y]\mapsto[x\otimes y]
$$
defines on $\H_{G\times H}(\cA)$ a module structure over the vertex algebra $\H_G(\C)\otimes\H_H(\C)$.}

\theorem{
Let $\cA$ be a half $O(\gs\gg\oplus\gs\gh)$ TVA which is separately $G$ and $H$ chiral free, let $\gg^*\oplus\gh^*\ra\cA$, $\zeta'\mapsto\theta_{\zeta'}$, be a $G\times H$ module homomorphism defining the $G$ and $H$ connection forms, and put $\Gamma_{\zeta'}=g\circ_0\theta_{\zeta'}$. Assume further that $G$ is semisimple, and that for $\zeta'\in\gg^*\oplus\gh^*$, $\zeta\in\gg\oplus\gh$, $\Gamma_{\zeta'}$ is $G\times H$ chiral horizontal and $\Gamma_{\zeta'}\circ_0\cA^0[0]=\iota_\zeta\circ_0\cA^0[0]=0$. Then the image of the class $\L\in\H_G(\C)$ in $\H_{G\times H}(\cA)$ is zero.}
\thmlab\GHCohomology

Note that since $g$ commutes with the vertex operators $\iota_\zeta$ in $\cA$ for $\zeta\in\gg\oplus\gh$, the assumption that $\Gamma_{\zeta'}=g\circ_0\theta_{\zeta'}\in\cA_{G\times H-hor}$ is equivalent to 
(cf. condition \dumbStar)
\eqn\dumb{
g\circ_0(\iota_\zeta\circ_0\theta_{\zeta'})=0,~~~\forall\zeta\in\gg\oplus\gh.
}
This holds a priori if $\zeta\in\gg$ and $\zeta'\in\gg^*$, for then we have $\iota_\zeta\circ_0\theta_{\zeta'}=\bra\zeta,\zeta'\ket$. The same holds for $\zeta\in\gh$, $\zeta'\in\gh'$. In any case, \dumb~and the assumption that $\Gamma_{\zeta'}\circ_0\cA^0[0]=\iota_\zeta\circ_0\cA^0[0]=0$
always hold in the geometrical setting $\cA=\cQ'(M)$. In fact, in this case we have $\cA^0[0]=C^\infty(M)$, giving $g\circ_0\cA^0[0]=0$, hence \dumb~ as well.

\corollary{Suppose, in addition to the preceding theorem, that $H=T$ is abelian. Then $\H_G(\cA_{T-bas})$ has no classes with positive weights.}
\proof
By Theorem \TReduction, we have an isomorphism $\H_G(\cA_{T-bas})\ra \H_{G\times T}(\cA)$ induced by the canonical inclusion map $\cW(\gg)\otimes\cA\hra\cW(\gt)\otimes\cW(\gg)\otimes\cA$. Under the natural map $\H_G(\C)\ra \H_{G\times T}(\cA)$ above, the image of the class $\L$ is $1\otimes\L\otimes 1$. Pulling back to $\H_G(\cA_{T-bas})$, it becomes $\L\otimes1\in\cW(\gg)\otimes\cA$. This coincides with $\kappa_G(\L)$, the image of $\L$ under the chiral Chern-Weil map $\kappa_G:\H_G(\C)\ra\H_G(\cA_{T-bas})$. Therefore, by the preceding theorem, 
this class must be zero. By Theorem \VirasoroVanishII, $\H_G(\cA_{T-bas})$ has no classes with positive weights. $\Box$

\corollary{Let $M$ be a $G\times H$ manifold in which the $G$ and $H$ actions on $M$ are free separately. Then the image of the class $\L\in\H_G(\C)$ in $\H_{G\times H}(\cQ'(M))$ is zero.}
\proof
Specialize Theorem \GHCohomology~to $\cA=\cQ'(M)$, which is separately $G$ and $H$ chiral free, since the $G$ and $H$ actions on $M$ are free separately. $\Box$

\corollary{Suppose, in addition to the preceding corollary, that $H=T$ is abelian. Then $\H_G(\cQ'(M/T))$ has no classes with positive weights.}
\proof
Specialize the first corollary to Theorem \GHCohomology~to the case $\cA=\cQ'(M)$. Then $\cA_{T-bas}\cong\cQ'(M/T)$ by Theorem \SheafGReduction. $\Box$

\corollary{For any semisimple group $G$, and a torus $T\subset G$ acting on $G$ by right multiplication, $\H_G(\cQ'(G/T))$ has no classes with positive weights.}

The theorem says in particular that $\kappa_G(\q)=0$ in $\H_G(\cQ'(G/T))$. Note that the $G$ action on $G/T$ in this case is never locally free. Yet the $G$ manifold $G/T$ has no weight one chiral equivariant cohomology. This is in sharp contrast to Theorem \AbelianFreeAction~for abelian $G$, where the vanishing of $\kappa_G(\q)$ guarantees that the $G$ action is locally free.

This example also illustrates an interesting distinction between the classical and the chiral theory. The former is a cohomology theory attached to pairs $(G,M)$ in a way that is (contravariantly) functorial in both $G$ and $M$. It is well-known that if the $G$ and $T$ actions are separately free on $M$, then we have the following natural isomorphisms of classical equivariant cohomology
$$
H_{G\times T}(M)\cong H_G(M/T)\cong H_T(M/G).
$$
Thus the classical theory cannot distinguish the three pairs $(G\times T,M),(G,M/T),(T,M/G)$.
But the chiral theory can sometimes distinguish at least the last two pairs. For $M=G$, and $T\subset G$ abelian, the preceding theorem shows that $\H_G(\cQ'(G/T))[1]=0$. But in \LLI, we have shown that $\H_T(\cQ'(G/G))=\H_T(\C)$ is nonzero in every positive weight. On the other hand, Corollary \TReductioncor~shows that the first two pairs lead to the same chiral equivariant cohomology, whether or not $G$ is abelian, and whether or not $G$ acts freely on $M$.

{\it Proof of Theorem \GHCohomology:} Fix orthonormal bases $\xi_i$ of $\gg$ and $\eta_j$ of $\gh$, and identify $\gg\equiv\gg^*$, $\gh\equiv\gh^*$, by using invariant forms.
Recall that the class $\L$ in $\H_G(\C)$ is represented by $d_{\cW(\gg)}(\Theta_\cS^{\xi_i}b^{\xi_i}+\beta^{\xi_i}\partial c^{\xi_i})$. Put
$$
z=(\Theta_\cS^{\xi_i}b^{\xi_i}+\beta^{\xi_i}\partial c^{\xi_i})\otimes 1\in C=\cW(\gg\oplus\gh)\otimes\cA.
$$ 
We shall adopt the convention that operators {\it operating} on either factor of $C$ will be regarded as operators acting on $C$. (In other words we shall omit tensoring such operators with 1.)
Suppose we can find an element $y\in C$ such that the following holds:
\item{(1)} $(\Theta_\cW^\xi+L_\xi)\circ_p \alpha=\beta^\xi\otimes 1~\delta_{p,1}$ for $\xi\in\gg$, $p\geq0$;
\item{(2)} $\alpha$ is $G$ chiral horizontal;
\item{(3)} $\alpha$ is $H$ chiral basic.

Put $d=d_{\cW(\gh)}+d_{\cW(\gg)}+d_\cA$. Then we claim that $z-d\alpha$ is $G\times H$ chiral horizontal. For $\zeta\in\gg\oplus\gh$, $p\geq0$, we have
$$
(b^\zeta+\iota_\zeta)\circ_p(z-d\alpha)
=\beta^{\zeta(\gg)}\otimes1~\delta_{p,1}-(\Theta_\cW^\zeta+L_\zeta)\circ_p \alpha+d[(b^\zeta+\iota_\zeta)\circ_p\alpha]
$$
where $\zeta(\gg)$ is the image of $\zeta$ under the projection $\gg\oplus\gh\ra\gg$. Consider the terms on the right side. The last term is zero by (2) and (3). When $\zeta\in\gg$, the first two terms add to zero by (1). When $\zeta\in\gh$, the first term is zero, and the second term is also zero by (3). This proves the claim.

If (1)-(3) hold, then we are done. For then the image of $\L$ in $\H_{G\times H}(\cA)$, which is $\L\otimes 1$, is represented by $dz=d(z-d\alpha)$. But since $z-d\alpha$ and $dz$ are both $G\times H$ chiral horizontal, it follows that $z-d\alpha$ is also $G\times H$ chiral invariant, by the identity $[d,b^\zeta+\iota_\zeta]=\Theta_\cW^\zeta+L_\zeta$. This means that $d(z-d\alpha)$ is $d$ exact in $(\cW(\gg\oplus\gh)\otimes\cA)_{G\times H-bas}$ and hence $\L\otimes 1$ is trivial.

We now proceed to construct $\alpha$ satisfying (1)-(3). First, in the simplest case $H=1$, Theorem \SemisimpleChiralFree~(for $G$ semisimple) already provides the answer: $\alpha=\beta^{\xi_i}\otimes\Gamma_{\xi_i}$. Starting from this, there will be a 4-step iteration in which we successively correct this to make (1)-(3) hold in general. Put
$$
y_0=\beta^{\xi_i}\otimes\Gamma_{\xi_i}.
$$
For $\alpha=y_0$, (1) holds because $y_0$ is $G$ invariant, $\Theta_\cW^\xi\circ_p$ kills the $\beta^{\xi_i}$ and $L_\xi\circ_p\Gamma_{\xi_i}=\bra\xi,\xi_i\ket~\delta_{p,1}$ for $\xi\in\gg$, $p>0$, by Lemma \GammaTransform; (2) also holds because the $\Gamma_{\xi_i}$ are $G$ chiral horizontal, by the same lemma; $y_0$ is also $T$ chiral horizontal because of the assumption $\Gamma_{\xi_i}\in\cA_{G\times H-hor}$. 

But for $\alpha=y_0$, (3) fails. Clearly $y_0$ is $H$ invariant, and it is killed by $(\Theta_\cW^\eta+L_\eta)\circ_p$ for $\eta\in\gh$, $p>1$, by weight consideration. But for $p=1$, this fails since
$$
(\Theta_\cW^\eta+L_\eta)\circ_1 y_0=\beta^{\xi_i}\otimes L_\eta\circ_1\Gamma_{\xi_i}.
$$
Note that $L_\eta\circ_1\Gamma_{\xi_i}=\iota_\eta\circ_0\theta_{\xi_i}\in\cA^0[0]$ by Lemma \GammaGTTransform. To cancel this
we put
$$
y_1=y_0+x_0,~~~x_0=-\beta^{\xi_i}\otimes:\Gamma_{\eta_j}(L_{\eta_j}\circ_1\Gamma_{\xi_i}):.
$$
Indeed, $x_0$ is manifestly $H$ invariant and is killed by $(\Theta_\cW^\eta+L_\eta)\circ_p$ for $\eta\in\gh$, $p>1$, by weight considerations. For $p=1$, we have
$$
(\Theta_\cW^\eta+L_\eta)\circ_1x_0=-\beta^{\xi_i}\otimes L_\eta\circ_1\Gamma_{\xi_i}
$$ 
by Lemma \GammaTransform, Lemma 2.9 of \LLI, and the assumption that $\Gamma_{[\eta,\eta_j]}\circ_0\cA[0]=0$. Thus $y_1=y_0+x_0$ is now $H$ chiral invariant. It is $G\times H$ chiral horizontal because $\Gamma_\zeta\in\cA_{G\times H-hor}$ and $\iota_\zeta\circ_0\cA^0[0]=0$, by assumptions. So for $\alpha=y_1$, (2) and (3) hold.

But for $\alpha=y_1$, (1) fails. In fact, we have seen that (1) holds for $y_0$, and $x_0$ is manifestly $G$ invariant and is killed by $(\Theta_\cW^\xi+L_\xi)\circ_p $ for $\xi\in\gg$, $p>1$, by weight consideration. For $p=1$, we have 
$$
(\Theta_\cW^\xi+L_\xi)\circ_1 x_0=-\beta^{\xi_i}\otimes (L_\xi\circ_1\Gamma_{\eta_j})(L_{\eta_j}\circ_1\Gamma_{\xi_i}).
$$
To cancel this we put
$$
y_2=y_1+x_1,~~~x_1=\Theta_\cE^{\xi_k}\beta^{\xi_i}\otimes(L_{\xi_k}\circ_1\Gamma_{\eta_j})(L_{\eta_j}\circ_1\Gamma_{\xi_i}).
$$
Indeed, $x_1$ is manifestly $G$ invariant and is killed by $(\Theta_\cW^\xi+L_\xi)\circ_p $ for $\xi\in\gg$, $p>1$, by weight consideration. For $p=1$, we have 
$$
(\Theta_\cW^\xi+L_\xi)\circ_1x_1=\beta^{\xi_i}\otimes (L_\xi\circ_1\Gamma_{\eta_j})(L_{\eta_j}\circ_1\Gamma_{\xi_i}).
$$
Here we have used the fact that $\Theta_\cW^\xi\circ_1\Theta_\cE^{\xi_k}=\Theta_\cE^\xi\circ_1\Theta_\cE^{\xi_k}=\bra\xi,\xi_k\ket$. So for $\alpha=y_2$, now (1) holds. We have seen that (3) holds for $y_1$. The term $x_1$ is also $H$ chiral basic, since all three vertex operators appearing in it, namely $\Theta_\cE^{\xi_k}$, $\beta^{\xi_i}$, $(L_{\xi_k}\circ_1\Gamma_{\eta_j})(L_{\eta_j}\circ_1\Gamma_{\xi_i})$, are $H$ chiral basic. Here we have used the fact that $G,H$ commute and that $\iota_\xi\circ_0\cA^0[0]=0$. So (3) holds for $\alpha=y_2$.

But for $\alpha=y_2$, (2) fails. In fact, we have seen that (2) holds for $y_1$, and $x_1$ is killed by $(b^\xi+\iota_\xi)\circ_p $ for $\xi\in\gg$, $p>0$, by weight consideration and because of $b^\xi(z)\Theta_\cE^{\xi_k}(w)\sim-b^{[\xi_k,\xi]}(w)~(z-w)^{-1}$. For $p=0$, we have 
$$
(b^\xi+\iota_\xi)\circ_0 x_1 =-b^{[\xi_k,\xi]}\beta^{\xi_i}\otimes(L_{\xi_k}\circ_1\Gamma_{\eta_j})(L_{\eta_j}\circ_1\Gamma_{\xi_i}).
$$
To cancel this we put
$$
y_3=y_2+x_2,~~~x_2=-b^{[\xi_k,\xi_l]}\beta^{\xi_i}\otimes\theta_{\xi_l}(L_{\xi_k}\circ_1\Gamma_{\eta_j})(L_{\eta_j}\circ_1\Gamma_{\xi_i}).
$$
Indeed, $x_2$ is killed by $(b^\xi+\iota_\xi)\circ_p $ for $\xi\in\gg$, $p>0$, by weight consideration and because $b^\xi$ commutes with $b^{[\xi_k,\xi_l]},\beta^{\xi_i}$. For $p=0$, we have 
$$
(b^\xi+\iota_\xi)\circ_0 x_2 =-b^{[\xi_k,\xi]}\beta^{\xi_i}\otimes (L_\xi\circ_1\Gamma_{\eta_j})(L_{\eta_j}\circ_1\Gamma_{\xi_i}).
$$
Here we have used the fact that $\iota_\xi\circ_0\theta_{\xi_l}=\bra\xi,\xi_l\ket$ and that $\iota_\xi\circ_0\cA^0[0]=0$. So for $\alpha=y_3$, now (2) holds. We have seen that (1) holds for $y_2$. The term $x_2$ is manifestly $G$ invariant, and the vertex operators appearing in it, namely $b^{[\xi_k,\xi_l]}$, $\beta^{\xi_i}$, $\theta_{\xi_l}(L_{\xi_k}\circ_1\Gamma_{\eta_j})(L_{\eta_j}\circ_1\Gamma_{\xi_i})$, are all killed by $(\Theta_\cW^\xi+L_\xi)\circ_p $ for $\xi\in\gg$, $p>0$, by weight consideration. So for $\alpha=y_3$, (1) holds.

But for $\alpha=y_3$, (3) fails. In fact, we have seen that (3) holds for $y_2$, and $x_2$ is manifestly $H$ chiral invariant. So $y_3$ is $H$ chiral invariant. The term $x_2$ is also killed by $(b^\eta+\iota_\eta)\circ_p $ for $\eta\in\gh$, $p>0$, by weight consideration and because $b^\eta$ commute with $b^{[\xi_k,\xi_l]},\beta^{\xi_i}$. For $p=0$, we have
$$
(b^\eta+\iota_\eta)\circ_0 x_2 =b^{[\xi_k,\xi_l]}\beta^{\xi_i}\otimes (L_\eta\circ_1\Gamma_{\xi_l})(L_\xi\circ_1\Gamma_{\eta_j})(L_{\eta_j}\circ_1\Gamma_{\xi_i}).
$$
Here we have used the fact that $\iota_\eta\circ_0\theta_{\xi_l}=L_\eta\circ_1\Gamma_{\xi_l}$, by Lemma \GammaGTTransform. 
To cancel this we put
$$
y_4=y_3+x_3,~~~ x_3=-c^{\eta_m}b^{[\xi_k,\xi_l]}\beta^{\xi_i}\otimes(L_{\eta_m}\circ_1\Gamma_{\xi_l})(L_{\xi_k}\circ_1\Gamma_{\eta_j})(L_{\eta_j}\circ_1\Gamma_{\xi_i}).
$$
(Note that all the vertex operators appearing in $x_3$ are pairwise commuting.)
Indeed, $x_3$ is killed by $(b^\eta+\iota_\eta)\circ_p $ for $\xi\in\gh$, $p>0$, by weight consideration and because $b^\eta$ commutes with $b^{[\xi_k,\xi_l]},\beta^{\xi_i}$. For $p=0$, we have 
$$
(b^\eta+\iota_\eta)\circ_0 x_3
=-b^{[\xi_k,\xi_l]}\beta^{\xi_i}\otimes(L_{\eta}\circ_1\Gamma_{\xi_l}) (L_\xi\circ_1\Gamma_{\eta_j})(L_{\eta_j}\circ_1\Gamma_{\xi_i}).
$$
Here we have used that $b^\eta\circ_0c^{\eta_m}=\bra\eta,\eta_m\ket$ and that $\iota_\eta\circ_0\cA^0[0]=0$. So $y_4$ is now $H$ chiral horizontal. We have seen that $y_3$ is $H$ chiral invariant. The term $x_3$ is also manifestly $H$ invariant and is killed by $(\Theta_\cW^\eta+L_\eta)\circ_p$ for $\eta\in\gh$, $p>0$, by weight consideration. So for $\alpha=y_4$, (3) holds, i.e. $y_4$ is $H$ chiral basic.

We have seen that (1) holds for $y_3$. The term $x_3$ is manifestly $G$ invariant, and the vertex operators appearing in it, namely $c^{\eta_m}$, $b^{[\xi_k,\xi_l]}$, $\beta^{\xi_i}$, $(L_{\eta_m}\circ_1\Gamma_{\xi_l})(L_{\xi_k}\circ_1\Gamma_{\eta_j})(L_{\eta_j}\circ_1\Gamma_{\xi_i})$, are all killed by $(\Theta_\cW^\xi+L_\xi)\circ_p $ for $\xi\in\gg$, $p>0$, by weight consideration. So for $\alpha=y_4$, (1) holds.

Finally, we have seen that (2) holds for $y_3$. The term $x_3$ is $G$ chiral horizontal,
since the vertex operators appearing in it are each $G$ chiral horizontal. Here we have used the fact that $\gg,\gh$ are orthogonal and that $\iota_\xi\circ_0\cA^0[0]$ for $\xi\in\gg$. So for $\alpha=y_4$, (2) holds.

In conclusion, for 
$$\eqalign{
\alpha &=\beta^{\xi_i}\otimes\Gamma_{\xi_i} - \beta^{\xi_i}\otimes:\Gamma_{\eta_j}(L_{\eta_j}\circ_1\Gamma_{\xi_i}):  +\Theta_\cE^{\xi_k}\beta^{\xi_i}\otimes(L_{\xi_k}\circ_1\Gamma_{\eta_j})(L_{\eta_j}\circ_1\Gamma_{\xi_i})\cr
&-c^{\eta_m}b^{[\xi_k,\xi_l]}\beta^{\xi_i}\otimes(L_{\eta_m}\circ_1
\Gamma_{\xi_l})(L_{\xi_k}\circ_1\Gamma_{\eta_j})(L_{\eta_j}\circ_1\Gamma_{\xi_i}),
}$$ 
we have shown that (1)-(3) hold. This completes the proof. $\Box$





\footatend\vfill\supereject\immediate\closeout\rfile\writestoppt
\baselineskip=14pt\centerline{{\bf References}}\bigskip{\frenchspacing%
\parindent=20pt\escapechar=` \input refs.tmp\vfill\eject}\nonfrenchspacing

\end

%% file: abbrev.tex

%
%
%
\def\unredoffs{} \def\redoffs{\voffset=-.31truein\hoffset=-.59truein}
\def\speclscape{\special{ps: landscape}}
%
%
%
%
\newbox\leftpage \newdimen\fullhsize \newdimen\hstitle \newdimen\hsbody
\tolerance=1000\hfuzz=2pt
\catcode`\@=11 
%
\ifx\answ\bigans\message{(This will come out unreduced.}
\magnification=1200\unredoffs\baselineskip=16pt plus 2pt minus 1pt
\hsbody=\hsize \hstitle=\hsize 
\else\message{(This will be reduced.} \let\l@r=L
\magnification=1000\baselineskip=16pt plus 2pt minus 1pt \vsize=7truein
\redoffs \hstitle=8truein\hsbody=4.75truein\fullhsize=10truein\hsize=\hsbody
\output={\ifnum\pageno=0 
  \shipout\vbox{\speclscape{\hsize\fullhsize\makeheadline}
    \hbox to \fullhsize{\hfill\pagebody\hfill}}\advancepageno
  \else
  \almostshipout{\leftline{\vbox{\pagebody\makefootline}}}\advancepageno
  \fi}
\def\almostshipout#1{\if L\l@r \count1=1 \message{[\the\count0.\the\count1]}
      \global\setbox\leftpage=#1 \global\let\l@r=R
 \else \count1=2
  \shipout\vbox{\speclscape{\hsize\fullhsize\makeheadline}
      \hbox to\fullhsize{\box\leftpage\hfil#1}}  \global\let\l@r=L\fi}
\fi
%
\newcount\yearltd\yearltd=\year\advance\yearltd by -1900

%

%
%

\def\draftmode{\message{ DRAFTMODE }\def\draftdate{{\rm preliminary draft:
\number\month/\number\day/\number\yearltd\ \ \hourmin}}%
\headline={\hfil\draftdate}\writelabels\baselineskip=20pt plus 2pt minus 2pt
 {\count255=\time\divide\count255 by 60 \xdef\hourmin{\number\count255}
  \multiply\count255 by-60\advance\count255 by\time
  \xdef\hourmin{\hourmin:\ifnum\count255<10 0\fi\the\count255}}}
\def\nolabels{\def\wrlabeL##1{}\def\eqlabeL##1{}\def\reflabeL##1{}}
\def\writelabels{\def\wrlabeL##1{\leavevmode\vadjust{\rlap{\smash%
{\line{{\escapechar=` \hfill\rlap{\sevenrm\hskip.03in\string##1}}}}}}}%
\def\eqlabeL##1{{\escapechar-1\rlap{\sevenrm\hskip.05in\string##1}}}%
\def\reflabeL##1{\noexpand\llap{\noexpand\sevenrm\string\string\string##1}}}
\nolabels
%
\global\newcount\secno \global\secno=0
\global\newcount\meqno \global\meqno=1
\def\newsec#1{\global\advance\secno by1\message{(\the\secno. #1)}
\global\subsecno=0\eqnres@t\noindent{\bf\the\secno. #1}
\writetoca{{\secsym} {#1}}\par\nobreak\medskip\nobreak}
\def\eqnres@t{\xdef\secsym{\the\secno.}\global\meqno=1\bigbreak\bigskip}
\def\sequentialequations{\def\eqnres@t{\bigbreak}}\xdef\secsym{}
\global\newcount\subsecno \global\subsecno=0
\def\subsec#1{\global\advance\subsecno by1\message{(\secsym\the\subsecno. #1)}
\ifnum\lastpenalty>9000\else\bigbreak\fi
\noindent{\it\secsym\the\subsecno. #1}\writetoca{\string\quad
{\secsym\the\subsecno.} {#1}}\par\nobreak\medskip\nobreak}
\def\appendix#1#2{\global\meqno=1\global\subsecno=0\xdef\secsym{\hbox{#1.}}
\bigbreak\bigskip\noindent{\bf Appendix #1. #2}\message{(#1. #2)}
\writetoca{Appendix {#1.} {#2}}\par\nobreak\medskip\nobreak}
%
%
\def\eqnn#1{\xdef #1{(\secsym\the\meqno)}\writedef{#1\leftbracket#1}%
\global\advance\meqno by1\wrlabeL#1}
\def\eqna#1{\xdef #1##1{\hbox{$(\secsym\the\meqno##1)$}}
\writedef{#1\numbersign1\leftbracket#1{\numbersign1}}%
\global\advance\meqno by1\wrlabeL{#1$\{\}$}}
\def\eqn#1#2{\xdef #1{(\secsym\the\meqno)}\writedef{#1\leftbracket#1}%
\global\advance\meqno by1$$#2\eqno#1\eqlabeL#1$$}
%
\newskip\footskip\footskip14pt plus 1pt minus 1pt 
\def\footnotefont{\ninepoint}\def\f@t#1{\footnotefont #1\@foot}
\def\f@@t{\baselineskip\footskip\bgroup\footnotefont\aftergroup\@foot\let\next}
\setbox\strutbox=\hbox{\vrule height9.5pt depth4.5pt width0pt}
\global\newcount\ftno \global\ftno=0
\def\foot{\global\advance\ftno by1\footnote{$^{\the\ftno}$}}
%
\newwrite\ftfile
\def\footend{\def\foot{\global\advance\ftno by1\chardef\wfile=\ftfile
$^{\the\ftno}$\ifnum\ftno=1\immediate\openout\ftfile=foots.tmp\fi%
\immediate\write\ftfile{\noexpand\smallskip%
\noexpand\item{f\the\ftno:\ }\pctsign}\findarg}%
\def\footatend{\vfill\eject\immediate\closeout\ftfile{\parindent=20pt
\centerline{\bf Footnotes}\nobreak\bigskip\input foots.tmp }}}
\def\footatend{}
%
%
\global\newcount\refno \global\refno=1
\newwrite\rfile
%
\def\ref{\nref}
\def\nref#1{\xdef#1{[\the\refno]}\writedef{#1\leftbracket#1}%
\ifnum\refno=1\immediate\openout\rfile=refs.tmp\fi
\global\advance\refno by1\chardef\wfile=\rfile\immediate
\write\rfile{\noexpand\item{#1\ }\reflabeL{#1\hskip.31in}\pctsign}\findarg}
\def\findarg#1#{\begingroup\obeylines\newlinechar=`\^^M\pass@rg}
{\obeylines\gdef\pass@rg#1{\writ@line\relax #1^^M\hbox{}^^M}%
\gdef\writ@line#1^^M{\expandafter\toks0\expandafter{\striprel@x #1}%
\edef\next{\the\toks0}\ifx\next\em@rk\let\next=\endgroup\else\ifx\next\empty%
\else\immediate\write\wfile{\the\toks0}\fi\let\next=\writ@line\fi\next\relax}}
\def\striprel@x#1{} \def\em@rk{\hbox{}}
\def\lref{\begingroup\obeylines\lr@f}
\def\lr@f#1#2{\gdef#1{\ref#1{#2}}\endgroup\unskip}

\def\addref#1{\immediate\write\rfile{\noexpand\item{}#1}} 
\def\startrefs#1{\immediate\openout\rfile=refs.tmp\refno=#1}
\def\refs#1{\count255=1[\r@fs #1{\hbox{}}]}
\def\r@fs#1{\ifx\und@fined#1\message{reflabel \string#1 is undefined.}%
\nref#1{need to supply reference \string#1.}\fi%
\vphantom{\hphantom{#1}}\edef\next{#1}\ifx\next\em@rk\def\next{}%
\else\ifx\next#1\ifodd\count255\relax\xref#1\count255=0\fi%
\else#1\count255=1\fi\let\next=\r@fs\fi\next}
%

%
\newwrite\ffile\global\newcount\figno \global\figno=1
\def\fig{fig.~\the\figno\nfig}
\def\nfig#1{\xdef#1{fig.~\the\figno}%
\writedef{#1\leftbracket fig.\noexpand~\the\figno}%
\ifnum\figno=1\immediate\openout\ffile=figs.tmp\fi\chardef\wfile=\ffile%
\immediate\write\ffile{\noexpand\medskip\noexpand\item{Fig.\ \the\figno. }
\reflabeL{#1\hskip.55in}\pctsign}\global\advance\figno by1\findarg}
\def\xfig{\expandafter\xf@g}\def\xf@g fig.\penalty\@M\ {}
\def\figs#1{figs.~\f@gs #1{\hbox{}}}
\def\f@gs#1{\edef\next{#1}\ifx\next\em@rk\def\next{}\else
\ifx\next#1\xfig #1\else#1\fi\let\next=\f@gs\fi\next}
\newwrite\lfile
{\escapechar-1\xdef\pctsign{\string\%}\xdef\leftbracket{\string\{}
\xdef\rightbracket{\string\}}\xdef\numbersign{\string\#}}

\def\writestop{\def\writestoppt{\immediate\write\lfile{\string\pageno%
\the\pageno\string\startrefs\leftbracket\the\refno\rightbracket%
\string\def\string\secsym\leftbracket\secsym\rightbracket%
\string\secno\the\secno\string\meqno\the\meqno}\immediate\closeout\lfile}}
\def\writestoppt{}\def\writedef#1{}
\def\seclab#1{\xdef #1{\the\secno}\writedef{#1\leftbracket#1}\wrlabeL{#1=#1}}
\def\subseclab#1{\xdef #1{\secsym\the\subsecno}%
\writedef{#1\leftbracket#1}\wrlabeL{#1=#1}}
\newwrite\tfile \def\writetoca#1{}
\def\leaderfill{\leaders\hbox to 1em{\hss.\hss}\hfill}
\def\writetoc{\immediate\openout\tfile=toc.tmp
   \def\writetoca##1{{\edef\next{\write\tfile{\noindent ##1
   \string\leaderfill {\noexpand\number\pageno} \par}}\next}}}
%
%
%

\catcode`\@=12 
%
\edef\tfontsize{\ifx\answ\bigans scaled\magstep3\else scaled\magstep4\fi}
\font\titlerm=cmr10 \tfontsize \font\titlerms=cmr7 \tfontsize
\font\titlermss=cmr5 \tfontsize \font\titlei=cmmi10 \tfontsize
\font\titleis=cmmi7 \tfontsize \font\titleiss=cmmi5 \tfontsize
\font\titlesy=cmsy10 \tfontsize \font\titlesys=cmsy7 \tfontsize
\font\titlesyss=cmsy5 \tfontsize \font\titleit=cmti10 \tfontsize
\skewchar\titlei='177 \skewchar\titleis='177 \skewchar\titleiss='177
\skewchar\titlesy='60 \skewchar\titlesys='60 \skewchar\titlesyss='60
\def\titlefont{\def\rm{\fam0\titlerm}
\textfont0=\titlerm \scriptfont0=\titlerms \scriptscriptfont0=\titlermss
\textfont1=\titlei \scriptfont1=\titleis \scriptscriptfont1=\titleiss
\textfont2=\titlesy \scriptfont2=\titlesys \scriptscriptfont2=\titlesyss
\textfont\itfam=\titleit \def\it{\fam\itfam\titleit}\rm}
 \ifx\answ\bigans\else scaled\magstep1\fi
\ifx\answ\bigans\else

 \font\absi=cmmi10 scaled\magstep1
\font\absis=cmmi7 scaled\magstep1 \font\absiss=cmmi5 scaled\magstep1
\font\abssy=cmsy10 scaled\magstep1 \font\abssys=cmsy7 scaled\magstep1
\font\abssyss=cmsy5 scaled\magstep1 
\skewchar\absi='177 \skewchar\absis='177 \skewchar\absiss='177
\skewchar\abssy='60 \skewchar\abssys='60 \skewchar\abssyss='60
\fi
\font\ninerm=cmr9 \font\sixrm=cmr6 \font\ninei=cmmi9 \font\sixi=cmmi6
\font\ninesy=cmsy9 \font\sixsy=cmsy6 \font\ninebf=cmbx9
\font\nineit=cmti9 \font\ninesl=cmsl9 \skewchar\ninei='177
\skewchar\sixi='177 \skewchar\ninesy='60 \skewchar\sixsy='60
\def\ninepoint{\def\rm{\fam0\ninerm}
\textfont0=\ninerm \scriptfont0=\sixrm \scriptscriptfont0=\fiverm
\textfont1=\ninei \scriptfont1=\sixi \scriptscriptfont1=\fivei
\textfont2=\ninesy \scriptfont2=\sixsy \scriptscriptfont2=\fivesy
\textfont\itfam=\ninei \def\it{\fam\itfam\nineit}\def\sl{\fam\slfam\ninesl}%
\textfont\bffam=\ninebf \def\bf{\fam\bffam\ninebf}\rm}
%
%

\hyphenation{anom-aly anom-alies coun-ter-term coun-ter-terms}
\def\inv{^{\raise.15ex\hbox{${\scriptscriptstyle -}$}\kern-.05em 1}}

\def\Dsl{\,\raise.15ex\hbox{/}\mkern-13.5mu D} 
\def\dsl{\raise.15ex\hbox{/}\kern-.57em\partial}

\def\lspace{\ifx\answ\bigans{}\else\qquad\fi}
\def\lbspace{\ifx\answ\bigans{}\else\hskip-.2in\fi} 
\def\boxeqn#1{\vcenter{\vbox{\hrule\hbox{\vrule\kern3pt\vbox{\kern3pt
    \hbox{${\displaystyle #1}$}\kern3pt}\kern3pt\vrule}\hrule}}}
\def\mbox#1#2{\vcenter{\hrule \hbox{\vrule height#2in
        \kern#1in \vrule} \hrule}}  
%

\def\darr#1{\raise1.5ex\hbox{$\leftrightarrow$}\mkern-16.5mu #1}

\def\half{{\textstyle{1\over2}}} 
\def\roughly#1{\raise.3ex\hbox{$#1$\kern-.75em\lower1ex\hbox{$\sim$}}}

%
%


\def\frac#1#2{{#1\over#2}}

\def\half{\frac12}

\def\journal#1&#2(#3){\unskip, #1~\bf #2 \rm(19#3) }
\def\andjournal#1&#2(#3){\sl #1~\bf #2 \rm (19#3) }

\def\bra#1{\left\langle #1\right|}
\def\ket#1{\left| #1\right\rangle}

\catcode`\@=11\def\slash#1{\mathord{\mathpalette\c@ncel{#1}}}
\overfullrule=0pt
\def\steepslash{\c@ncel}
\def\frac#1#2{{#1\over #2}}

\def\:{\!:\!}
\def\inbar{\,\vrule height1.5ex width.4pt depth0pt}
\def\IQ{\relax\,\hbox{$\inbar\kern-.3em{\rm Q}$}}
\def\IB{\relax{\rm I\kern-.18em B}}
\def\IC{\relax\hbox{$\inbar\kern-.3em{\rm C}$}}
\def\IP{\relax{\rm I\kern-.18em P}}
\def\IR{\relax{\rm I\kern-.18em R}}
\def\ZZ{\relax\ifmmode\mathchoice
{\hbox{Z\kern-.4em Z}}{\hbox{Z\kern-.4em Z}}
{\lower.9pt\hbox{Z\kern-.4em Z}}
{\lower1.2pt\hbox{Z\kern-.4em Z}}\else{Z\kern-.4em Z}\fi}

\catcode`\@=12

\def\npb#1(#2)#3{{ Nucl. Phys. }{B#1} (#2) #3}
\def\plb#1(#2)#3{{ Phys. Lett. }{#1B} (#2) #3}
\def\pla#1(#2)#3{{ Phys. Lett. }{#1A} (#2) #3}
\def\prl#1(#2)#3{{ Phys. Rev. Lett. }{#1} (#2) #3}
\def\mpla#1(#2)#3{{ Mod. Phys. Lett. }{A#1} (#2) #3}
\def\ijmpa#1(#2)#3{{ Int. J. Mod. Phys. }{A#1} (#2) #3}
\def\cmp#1(#2)#3{{ Comm. Math. Phys. }{#1} (#2) #3}
\def\cqg#1(#2)#3{{ Class. Quantum Grav. }{#1} (#2) #3}
\def\jmp#1(#2)#3{{ J. Math. Phys. }{#1} (#2) #3}
\def\anp#1(#2)#3{{ Ann. Phys. }{#1} (#2) #3}
\def\prd#1(#2)#3{{ Phys. Rev. } {D{#1}} (#2) #3}
\def\ptp#1(#2)#3{{ Progr. Theor. Phys. }{#1} (#2) #3}
\def\aom#1(#2)#3{{ Ann. Math. }{#1} (#2) #3}

\def\bs{\bigskip}
\def\bu{\bs\noindent $\bullet$}

\def\br{\buildrel}
\def\bra{\langle}
\def\ket{\rangle}

\def\C{{\bf C}}

\def\H{{\bf H}}

\def\L{{\bf L}}

\def\R{{\bf R}}

\def\Z{{\bf Z}}
\def\cA{{\cal A}}
\def\cB{{\cal B}}
\def\cC{{\cal C}}
\def\cD{{\cal D}}
\def\cE{{\cal E}}
\def\cF{{\cal F}}

\def\cJ{{\cal J}}
\def\cK{{\cal K}}

\def\cP{{\cal P}}
\def\cQ{{\cal Q}}

\def\cS{{\cal S}}

\def\cV{{\cal V}}
\def\cW{{\cal W}}

\input amssym

\def\gg{{\goth g}}
\def\gh{{\goth h}}

\def\gs{{\goth s}}
\def\gt{{\goth t}}

\def\gX{{\goth X}}

\def\cicy#1(#2|#3)#4{\left(\matrix{#2}\right|\!\!
                     \left|\matrix{#3}\right)^{{#4}}_{#1}}

\def\emptyset{\varnothing}

\def\lra{\longrightarrow}

\def\hra{\hookrightarrow}
\def\hla{\hookleftarrow}
\def\ra{\rightarrow}
\def\la{\leftarrow}

\def\bs{\bigskip}

\def\Box{{\,\lower0.9pt\vbox{\hrule
\hbox{\vrule height 0.2 cm \hskip 0.2 cm
\vrule height 0.2 cm}\hrule}\,}}

\global\newcount\thmno \global\thmno=0
\def\definition#1{\global\advance\thmno by1
\bigskip\noindent{\bf Definition \secsym\the\thmno. }{\it #1}
\par\nobreak\medskip\nobreak}
\def\question#1{\global\advance\thmno by1
\bigskip\noindent{\bf Question \secsym\the\thmno. }{\it #1}
\par\nobreak\medskip\nobreak}
\def\theorem#1{\global\advance\thmno by1
\bigskip\noindent{\bf Theorem \secsym\the\thmno. }{\it #1}
\par\nobreak\medskip\nobreak}
\def\proposition#1{\global\advance\thmno by1
\bigskip\noindent{\bf Proposition \secsym\the\thmno. }{\it #1}
\par\nobreak\medskip\nobreak}
\def\corollary#1{\global\advance\thmno by1
\bigskip\noindent{\bf Corollary \secsym\the\thmno. }{\it #1}
\par\nobreak\medskip\nobreak}
\def\lemma#1{\global\advance\thmno by1
\bigskip\noindent{\bf Lemma \secsym\the\thmno. }{\it #1}
\par\nobreak\medskip\nobreak}
\def\conjecture#1{\global\advance\thmno by1
\bigskip\noindent{\bf Conjecture \secsym\the\thmno. }{\it #1}
\par\nobreak\medskip\nobreak}
\def\exercise#1{\global\advance\thmno by1
\bigskip\noindent{\bf Exercise \secsym\the\thmno. }{\it #1}
\par\nobreak\medskip\nobreak}
\def\remark#1{\global\advance\thmno by1
\bigskip\noindent{\bf Remark \secsym\the\thmno. }{\it #1}
\par\nobreak\medskip\nobreak}
\def\problem#1{\global\advance\thmno by1
\bigskip\noindent{\bf Problem \secsym\the\thmno. }{\it #1}
\par\nobreak\medskip\nobreak}
\def\others#1#2{\global\advance\thmno by1
\bigskip\noindent{\bf #1 \secsym\the\thmno. }{\it #2}
\par\nobreak\medskip\nobreak}
\def\proof{\noindent Proof: }

\def\thmlab#1{\xdef #1{\secsym\the\thmno}\writedef{#1\leftbracket#1}\wrlabeL{#1=#1}}
%
%
\def\newsec#1{\global\advance\secno by1\message{(\the\secno. #1)}
\global\subsecno=0\thmno=0\eqnres@t\noindent{\bf\the\secno. #1}
\writetoca{{\secsym} {#1}}\par\nobreak\medskip\nobreak}
\def\eqnres@t{\xdef\secsym{\the\secno.}\global\meqno=1\bigbreak\bigskip}
\def\sequentialequations{\def\eqnres@t{\bigbreak}}\xdef\secsym{}
%

%
\newcount{\exnum}
\def\prob{\advance\exnum by 1
\bigskip\item{\the\exnum.}\ }
\newcount{\exnum}
\def\next{\advance\exnum by 1
\bigskip\noindent{\the\exnum.}\ }
\def\np{\vfill\eject}